\documentclass[11pt,leqno]{article}

\usepackage{amsfonts,latexsym, amsmath, amssymb, amsthm, hyperref}

\usepackage{gfsartemisia}
\usepackage[T1]{fontenc}

\usepackage{color}

\newtheorem{theorem}{Theorem}[section]
\newtheorem{lemma}[theorem]{Lemma}

\newtheorem{proposition}[theorem]{Proposition}
\newtheorem{corollary}[theorem]{Corollary}
\newtheorem{remarks}[theorem]{Remarks}
\newtheorem{remark}[theorem]{Remark}

\theoremstyle{definition}

\newtheorem{note}[theorem]{Note}

\makeatletter

\numberwithin{equation}{section}

\makeatother

\DeclareMathOperator{\signum}{sgn}

\newcommand{\sgn}[1]{\signum(#1)}
\newcommand{\prend}{$\hfill \quad \Box$}

\newcommand\blfootnote[1]{%
  \begingroup
  \renewcommand\thefootnote{}\footnote{#1}%
  \addtocounter{footnote}{-1}%
  \endgroup
}

\begin{document}

\small

\title{Random version of Dvoretzky's theorem in $\ell_p^n$}

\author{Grigoris Paouris\thanks{Supported by the A. Sloan foundation, BSF grant 2010288 and the NSF
CAREER-1151711 grant;}, \, Petros Valettas\thanks{Supported in part by the NSF Grant DMS-1612936;} \, and \, Joel Zinn\thanks{Partially supported by NSF grants DMS-1208962 and DMS-1205781.}}

\maketitle

\begin{abstract}\footnotesize
We study the dependence on $\varepsilon$ in the critical dimension $k(n,p,\varepsilon)$ for which one can find random sections 
of the $\ell_p^n$-ball which are $(1+\varepsilon)$-spherical. We give lower (and upper) estimates for $k(n,p,\varepsilon)$
for all eligible values $p$ and $\varepsilon$ as $n\to \infty$, which agree with the sharp estimates for the extreme values $p=1$ and $p=\infty$. Toward this end, we provide tight bounds for the Gaussian concentration of the $\ell_p$-norm.
\end{abstract}



\section{Introduction}

The fundamental theorem of Dvoretzky from \cite{Dvo} in geometric language states that every centrally 
symmetric convex body on $\mathbb R^n$ has a central section of large dimension which is almost spherical. 
The optimal form of the theorem, which was proved by Milman 
in \cite{Mil}, reads as follows. For any $\varepsilon \in (0,1)$ there exists $\eta=\eta(\varepsilon)>0$ with
the following property: for every $n$-dimensional symmetric convex body $A$ there exist a linear image $A_1$ of $A$ 
and $k$-dimensional subspace $F$ with $k \geq \eta(\varepsilon) \log n$ such that
\begin{align*}
(1-\varepsilon) B_F \subseteq A_1\cap F \subseteq (1+\varepsilon)B_F,
\end{align*} where $B_F$ denotes the Euclidean ball in $F$. The example of the cube $A= B_\infty^n$ shows that this result is best 
possible with respect to $n$ (see \cite{Sch3} for the details). 
The approach of \cite{Mil} is
probabilistic in nature and shows that most of the $k$-dimensional sections 
are $(1+\varepsilon)$-spherical (or Euclidean). Here ``most" means with overwhelming probability in terms of the Haar probability measure $\nu_{n,k}$ 
on the Grassmann manifold $G_{n,k}$. More precisely, given a centrally symmetric convex body $A$ on $\mathbb R^n$ and $\varepsilon \in (0,1)$
the random $k$-dimensional subspace $F$ satisfies:
\begin{align*}
\frac{1-\varepsilon}{M} B_F \subseteq A\cap F \subseteq \frac{1+\varepsilon}{M}B_F
\end{align*} with probability greater than $1-e^{-k}$ as long as $k\leq c(\varepsilon)k(A)$. Here $c(\varepsilon)$ 
stands for the function of $\varepsilon$ in the probabilistic formulation and $k(A)$ is usually referred to the ``critical dimension"
of the body $A$. The latter can be computed in terms of the global parameters $M=M(A)=\int_{S^{n-1}} \|\theta\|_A \, d\sigma(\theta)$ 
and $b=b(A)=\max_{\theta \in S^{n-1}} \|\theta\|_A$; that is $k(A) \simeq n(M/b)^2$.  Recall that $1/b$ is the radius of the maximal centered
inscribed ball in $A$. Next, one may select a good position of the body
$A$ for which the $k(A)$ is large enough with respect to $n$ (see \cite{MS} for further details).

It has been proved in \cite{MS2} that this 
formulation is optimal with respect to the dimension $k(A)$ in the following sense: the maximal dimension $m$ for which the random
$m$-dimensional sections are $4$-Euclidean with probability greater than $\frac{n}{n+m}$ is less than $Ck(A)$ for some absolute constant $C>0$, i.e. $m\lesssim k(A)$.\footnote{For any two quantities $\Gamma,\Delta$ depending on $n,p$, etc. we write $\Gamma\lesssim \Delta$ if there exists numerical constant 
$C>0$ - independent of everything - such that $\Gamma \leq C \Delta$. We write $\Gamma \gtrsim \Delta$ if $\Delta \lesssim \Gamma$ and 
$\Gamma \simeq \Delta$ if $\Gamma \lesssim \Delta$ and 
$\Delta \lesssim \Gamma$. Accordingly we write $\Gamma \simeq_p \Delta$ if the constants involved are depending only on $p$.} (Here and everywhere else $C,c,C_1, c_1, \ldots$ stand for positive absolute constants whose values may change from line to line).

The proof in \cite{Mil} provides the lower bound $c(\varepsilon) \geq c\varepsilon^2 /\log \frac{1}{\varepsilon}$ 
and this is improved to $c(\varepsilon) \geq c\varepsilon^2$ 
by Gordon in \cite{Go} and an alternative approach is given by Schechtman in \cite{Sch1}. 
This dependence is known to be optimal. The recent works of Schechtman in \cite{Sch2} and Tikhomirov in \cite{Tik} established 
that the dependence on $\varepsilon$ in the randomized Dvoretzky for $B_\infty^n$ is of the exact order $\varepsilon/ \log \frac{1}{\varepsilon}$. 

As far as the dependence on $\varepsilon$ in the existential version of Dvoretzky's theorem is concerned, 
Schechtman proved in \cite{Sch2} that one can always $(1+\varepsilon)$-embed $\ell_2^k$ in any
$n$-dimensional normed space $E$ with $k(E, \varepsilon) \geq c\varepsilon \log n / (\log\frac{1}{\varepsilon})^2 $. Tikhomirov in \cite{Tik2} 
proved that for 1-symmetric spaces $E$ we may have $k(E,\varepsilon) \geq c \log n/ \log \frac{1}{\varepsilon}$ complementing the previously known result due to Bourgain and Lindenstrauss from \cite{BL}. Recall that a normed space $(\mathbb R^n,\|\cdot\|)$ is said to be 1-symmetric if the norm satisfies
$\|\sum_i \varepsilon_i a_i e_{\pi(i)}\|= \|\sum_i a_i e_i\|$ for all scalars $(a_i)$, for all choices of signs $\varepsilon_i=\pm 1$ and for any permutation $\pi$, where $(e_i)$ is the 
standard basis in $\mathbb R^n$. Tikhomirov's result was subsequently extended by
Fresen in \cite{Fres} for permutation invariant spaces with uniformly bounded basis constant. In this note we will not deal with the existential form of Dvoretzky's theorem. 
Related results for $\ell_p$ spaces are presented in \cite{Ko}. For more detailed information on the subject, explicit statements and 
historical remarks the reader is referred to the recent monograph \cite{AGM}.

Our goal here is to study the random version for the spaces $\ell_p^n$ and to give bounds on the dimension $k(n,p,\varepsilon)\equiv k(\ell_p^n, \varepsilon)$ for which the $k$-dimensional random section of $B_p^n$ is $(1+\varepsilon)$-Euclidean with high probability on $G_{n,k}$. 
These bounds are continuous with respect to $p$ and coincide with the known bounds in the extreme cases $p=1$ and $p=\infty$. 
To this end we first study the concentration phenomenon for the $\ell_p$ norms and we prove the following result:

\begin{theorem} \label{thm: 1.1} For all sufficiently large $n$ and for any $1\leq p\leq \infty$ one has:
\begin{align*}
P\left( \big | \|X\|_p-\mathbb E\|X\|_p \big| > \varepsilon \mathbb E\|X\|_p \right) \leq C_1\exp(-c_1\beta(n,p,\varepsilon)), \quad 0<\varepsilon<1,
\end{align*} where $X$ is standard $n$-dimensional Gaussian vector and $C_1,c_1>0$ are absolute constants. 
The function $\beta(n,p,\varepsilon)$ is defined as follows:
\begin{align*}
\beta(n,p,\varepsilon) =  \left\{
\begin{array}{lll}
\varepsilon^2n, & 1\leq p\leq 2 \\
 \max\left\{  \min \left\{ p^2 2^{-p}\varepsilon^2 n, (\varepsilon n)^{2/p} \right\}, \varepsilon pn^{2/p} \right\}, & 2<p\leq c_0\log n \\
\varepsilon pn^{2/p}, & p> c_0 \log n 
\end{array}\right. ,
\end{align*}
where $0< c_0 < 1$ is suitable absolute constant. Furthermore, for $p\leq c_0\log n$ we have:
\begin{align*}
P\left( \big | \|X\|_p-\mathbb E\|X\|_p \big| > \varepsilon \mathbb E\|X\|_p \right) \leq \frac{C_1}{1+p^2 2^{-p} \varepsilon^2 n},
\end{align*} for all $\varepsilon> 0$.
\end{theorem}

The bound we retrieve in the case of fixed $p$ is not new. The corresponding estimates have been studied by Naor \cite{Naor} in an
even more general probabilistic context. Also, for $p=\infty$ we recover the same bound proved by Schechtman in \cite{Sch2}. Therefore,
the above concentration result interpolates between the sharp concentration estimates for fixed $1\leq p<\infty$ and $p=\infty$ and
is derived in a unified way. However, our methods are different from the techniques used in \cite{Naor} and \cite{Sch2} and utilize 
Gaussian functional inequalities. Actually, following the same ideas as in \cite{Sch1} we will prove a distributional inequality for Gaussian 
random matrices similar to the concentration inequality
described above. Using this inequality and a chaining argument we prove the second main result which is the critical 
dimension $k(n,p,\varepsilon)$ in the randomized Dvoretzky for the $B_p^n$ balls.

\begin{theorem} \label{thm: 1.2} For all large enough $n$, for any $1\leq p\leq \infty$ and for any $0 < \varepsilon <1 $ the random $k$-dimensional section 
of $B_p^n$ with dimension $k\leq k(n,p,\varepsilon)$ is $(1+\varepsilon)$-Euclidean with probability greater than 
$1-C\exp(-c k(n,p,\varepsilon) )$, where $k(n,p,\varepsilon)$ is defined as:
\begin{itemize}
\item [i.] If $1\leq p<2$, then 
\begin{align*}
k(n,p,\varepsilon) \gtrsim  \varepsilon^2n.
\end{align*}

\item [ii.] If $2<p<c_0 \log n$, then
\begin{align}
k(n,p,\varepsilon) \gtrsim \left\{ 
\begin{array}{ccc}
(Cp)^{-p}  \varepsilon^2n , & 0< \varepsilon \leq (Cp)^{p/2} n^{-\frac{p-2}{2(p-1)}} \\
p^{-1} \varepsilon^{2/p} n^{2/p}, & (Cp)^{p/2} n^{-\frac{p-2}{2(p-1)}} <\varepsilon \leq 1/p \\
\varepsilon pn^{2/p}/\log \frac{1}{\varepsilon}, & 1/p < \varepsilon <1
\end{array}
\right. .
\end{align} Furthermore for $p< c_0\log n$ we have:
\begin{align*}
k(n,p,\varepsilon) \gtrsim \log n/ \log\frac{1}{\varepsilon}.
\end{align*}
\item [iii.] If $p\geq c_0 \log n$, then 
\begin{align*}
k(n,p, \varepsilon) \gtrsim  \varepsilon \log n/ \log\frac{1}{\varepsilon} .
\end{align*}
\end{itemize}
where $C,c, c_0>0$ are absolute constants. 
\end{theorem}

As one observes the dependence on $\varepsilon$ in $1\leq p\leq 2$ is $\varepsilon^2$ as predicted by V. Milman's proof (and its improvement by \cite{Gor} and \cite{Sch1}). However, for $p>2$ the
dependence on $\varepsilon$ is much better than $\varepsilon^2$ for all values of $p$. This permits us to find sections of $B_p^n$ of 
polynomial dimension which are closer to the Euclidean ball than previously obtained. Observe that Theorem \ref{thm: 1.2} retrieves the 
right dependence on $c(\varepsilon)$ at $p=1$ (actually when $p$ is fixed) and at $p=\infty$.  

The rest of the paper is organized as follows: In Section 2 we fix the notation, we give the required background material and 
we include some basic probabilistic inequalities. Gaussian  functional inequalities as logarithmic Sobolev inequality,
Talagrand's $L_1-L_2$ inequality and Pisier's Gaussian inequality are also included.
 
Before the proof of Theorem \ref{thm: 1.1} we prefer to deal with an easier problem first; the problem of determining
the right order of the Gaussian variance of the $\ell_p$ norm. We study this question in Section 3. This is a warm-up for the concentration
result we will investigate in Section 4. The main techniques that we will use, as well as the main problems we have to resolve, will be
apparent already in Section 3. This estimate will be used to obtain the dependence $\log n /\log \frac{1}{\varepsilon}$ 
for $p\leq c_0\log n$, but still proportional to $\log n$ in Theorem \ref{thm: 1.2}.

In Section 4 we present the proof of Theorem \ref{thm: 1.1}. Moreover, efforts have been made to provide lower estimates for the
probability described in Theorem \ref{thm: 1.1} (see also the Appendix by Tikhomirov). 

In Section 5 we prove Theorem \ref{thm: 1.2} and we show that in several cases the result is best possible up to constants. 

We conclude in Section 6 with further remarks and open questions.

\section{Notation and background material}

We work in $\mathbb R^n$ equipped with the standard inner product $\langle x, y\rangle =\sum_{i=1}^n x_iy_i$ for 
$x=(x_1,\ldots,x_n)$ and $y=(y_1,\ldots, y_n) $ in $\mathbb R^n$. The $\ell_p$-norm in $\mathbb R^n$ ($1\leq p<\infty$) is defined as:
\begin{align*}
\|x\|_{\ell_p^n} \equiv \|x\|_p := \left( \sum_{i=1}^n |x_i|^p\right)^{1/p} , \; x=(x_1,\ldots,x_n)
\end{align*} and for $p=\infty$ as:
\begin{align*}
\|x\|_{\ell_\infty^n} \equiv \|x\|_\infty := \max_{1\leq i\leq n} |x_i| , \, x=(x_1,\ldots, x_n).
\end{align*} The Euclidean sphere is defined as: $S^{n-1}= \{x\in \mathbb R^n : \|x\|_2=1 \}$.
The normed space $(\mathbb R^n,  \|\cdot\|_p)$ is denoted by $\ell_p^n$, for $1\leq p\leq \infty$ and its unit ball by $B_p^n$, 
i.e. $B_p^n=\{x\in \mathbb R^n : \|x\|_p\leq 1\}$.
For $1\leq p < q \leq \infty$ we have:
\begin{align}\label{eq: Holder - p-norms}
\|x\|_q \leq \|x\|_p \leq n^{1/p-1/q} \|x\|_q,
\end{align} for all $x\in \mathbb R^n$. We write $\|\cdot\|$ for an arbitrary norm on $\mathbb R^n$ and $\|\cdot\|_A$ if the
norm is induced by the centrally symmetric convex body $A$ on $\mathbb R^n$. For any subspace $F$ of $\mathbb R^n$ we write:
$S_F: = S^{n-1}\cap F$  and $B_F:= B_2^n \cap F$. 

The random variables in some probability space $(\Omega, \mathcal A, P)$ are denoted by $\xi , \eta, \ldots$ while the random vectors by $X=(X_1,\ldots, X_n)$ or simply $X, Y, Z, \ldots$. The random vectors under consideration are going to be Gaussian unless it is 
stated otherwise. If $\mu$ is a probability measure we write $\mathbb E_\mu$ and $\rm Var_\mu$ for the 
expectation and the variance respectively with respect to $\mu$. If the measure is prescribed the subscript is omitted. We shall make frequent use of the Paley-Zygmund inequality (for a proof see \cite{BLM}):

\begin{lemma} \label{lem: PZ-ineq}
Let $\xi$ be a non-negative random variable defined on some probability space $(\Omega, \mathcal A, P)$ 
with $\xi\in L_2(\Omega, \mathcal A, P)$. Then, 
\begin{align*}
P \left( \xi \geq t \mathbb E\xi \right) \geq (1-t)^2 \frac{(\mathbb E\xi)^2}{\mathbb E\xi^2},
\end{align*} for all $0<t<1$.
\end{lemma} 

\medskip

Also the multivariate version of Chebyshev's association inequality due to Harris will be useful:

\begin{proposition} \label{prop:Harris}  Let $Z=(\zeta_1,\ldots, \zeta_k)$ where $\zeta_1,\ldots,\zeta_k$ are i.i.d. 
random variables taking values almost surely in $A \subseteq \mathbb R$. 
If $F,G: A^k \subseteq \mathbb R^k \to \mathbb R$ are coordinatewise non-decreasing\footnote{A real valued function $H$ defined
on $U \subseteq \mathbb R^k$ is said to be {\it coordinatewise non-decreasing} if it is non-decreasing in each variable while keeping all the 
other variables fixed at any value.} functions, then we have:
\begin{align*}
\mathbb E [F(Z) G(Z) ] \geq \mathbb E [ F(Z) ] \mathbb E [ G(Z) ].
\end{align*}
\end{proposition} 

Harris' inequality can be derived from consecutive applications of Chebyshev's association inequality and conditioning. 
For a detailed proof we refer the reader to \cite{BLM}. For some measure space $( \Omega, \mathcal {E}, \mu)$ we write
\begin{align*}
\|f\|_{L_p(\mu)} := \left( \int_\Omega |f|^p \, d\mu \right)^{1/p}, \; 1\leq p<\infty.
\end{align*} for any measurable function $f:\Omega \to \mathbb R$. If $\mu$ is Borel probability measure on $\mathbb R^n$ and 
$K$ is a centrally symmetric convex body on $\mathbb R^n$ we also use the notation
\begin{align*}
I_r(\mu, K) := \left( \int_{\mathbb R^n} \|x\|_K^r \, d\mu(x) \right)^{1/r}, \, -n<r \neq 0,
\end{align*} and for $r=0$
\begin{align*}
I_0(\mu, K):= \exp\left( \int_{\mathbb R^n} \log \|x\|_K \, d\mu(x) \right).
\end{align*} If $\sigma$ is the (unique) probability measure on $S^{n-1}$ which is invariant under orthogonal transformations and 
$A$ is centrally symmetric convex body on $\mathbb R^n$, then we write:
\begin{align}
M_q(A) := \left( \int_{S^{n-1}} \|\theta\|_A^q \, d\sigma(\theta) \right)^{1/q}, \quad q\neq 0.
\end{align} For $q=1$ we simply write $M(A)=M_1(A)$.

For the random version of Dvoretzky's theorem recall V. Milman's formulation from \cite{Mil} (see also \cite{MS} or \cite{AGM}) and
see \cite{Gor} and \cite{Sch1} for the dependence on $\varepsilon$:

\begin{theorem} \label{thm: VMil}
 Let $A$ be a centrally symmetric convex body on $\mathbb R^n$. Define the critical dimension $k(A)$ of $A$ as follows:
 \begin{align} \label{eq:dvo-num}
 k(A) = \frac{ \mathbb E\|Z\|_A^2 }{b^2(A)} \simeq n\left( \frac{M(A)}{b(A)}\right)^2,
 \end{align} where $b(A)$ is the Lipschitz constant of the map $x\mapsto \|x\|_A$, i.e. $b=\max_{\theta\in S^{n-1}} \|\theta\|_A$ and $Z$
 is a standard Gaussian $n$-dimensional random vector. Then,
 the random $k$-dimensional subspace $F$ of $ (\mathbb R^n, \|\cdot\|_A)$ satisfies:
 \begin{align*}
 \frac{1}{(1+\varepsilon)M}B_F  \subseteq A\cap F \subseteq \frac{1}{(1-\varepsilon)M} B_F
 \end{align*} with probability greater than $1-e^{-ck}$ provided that $k\leq k(A,\varepsilon)$, where $k(A,\varepsilon) \simeq \varepsilon^2 k(A)$ and $M\equiv M(A)$.
\end{theorem}  Here the probability is considered with respect to the Haar probability measure $\nu_{n,k}$ on the Grassmann manifold
$G_{n,k}$, which is invariant under the orthogonal group action.

With some abuse of terminology for a subspace $F$ of a normed space $(\mathbb R^n, \|\cdot\|)$ (or equivalently for a section 
$A\cap F$ of a centrally symmetric convex body $A$ on $\mathbb R^n$) we say that is $(1+\varepsilon)$-{\it spherical} (or {\it Euclidean}) if:
\begin{align*}
\max_{\theta \in S_F} \|\theta\| / \min_{\theta \in S_F} \|\theta\| < 1+\varepsilon \quad {\rm or} 
\quad \max_{z\in S_F} \|z\|_A / \min_{z\in S_F} \|z\|_A <1+\varepsilon.
\end{align*} Thus, the previous theorem states that the random $k$-dimensional subspace of $(\mathbb R^n, \|\cdot\|_A)$ is $\frac{1+\varepsilon}{1-\varepsilon}$-spherical 
with probability greater than $1-e^{-ck}$ as long as $k\leq \varepsilon^2 k(A)$.
In the next paragraph we provide asymptotic estimates for $k_{p,n}:=k(\ell_p^n) \equiv k(B_p^n)$ in terms of 
$n$ and $p$.

\subsection{Gaussian random variables} 

If $g$ is a standard Gaussian random variable we set $\sigma_p^p:=\mathbb E|g|^p$
for every $p>0$. The next asymptotic estimate follows easily by Stirling's formula:
\begin{align} \label{eq: 2.4}
\sigma_p^p = \mathbb E|g|^p =\frac{2^{p/2}}{\sqrt{\pi}} \Gamma \left(\frac{p+1}{2} \right) 
\sim \sqrt{2} \left( \frac{p}{e} \right)^{p/2} , \quad p\to\infty.\end{align}

The $n$-dimensional standard Gaussian measure with density $(2\pi)^{-n/2} e^{-\|x\|_2^2/2}$ is denoted by $\gamma_n$. In the 
next Proposition, the asymptotic estimate \eqref{eq:mean-ell-p} is a special case of a more general result from \cite{SZ}.

\begin{proposition}\label{prop:mean-ell-p}
Let $1\leq p\leq \infty$ and let $Z$ be distributed according to $\gamma_n$. Then, we have:
\begin{align} \label{eq:mean-ell-p}
\mathbb E\|Z\|_p = \int_{\mathbb R^n} \|x\|_p \, d\gamma_n(x) \simeq  \left\{ \begin{array}{ll} 
n^{1/p}\sqrt{p}, & p< \log n\\
\sqrt{\log n}, & p\geq \log n 
\end{array} \right. .
\end{align} Therefore, for the critical dimension of $B_p^n$, we have:
\begin{align*}
k_{p,n}=k(B_p^n) \simeq \left\{ \begin{array}{lll}
n & 1\leq p\leq 2 \\
pn^{2/p} & 2\leq p \leq \log n \\
\log n & p\geq \log n
\end{array} \right. .
\end{align*}
\end{proposition}

We shall need Gordon's lemma for Mill's ratio from \cite{Go}:

\begin{lemma} \label{lem:Gordon-ineq} For any $a>0$ we have:
\begin{align} \label{eq:gordon-1}\frac{a}{1+a^2}\leq e^{a^2/2}\int_a^\infty e^{-t^2/2}\, dt\leq \frac{1}{a}. 
\end{align} Equivalently, we have:
\begin{align} \label{eq:gordon-2}
1\leq \frac{\phi(a)}{a(1-\Phi(a) )} \leq 1+\frac{1}{a^2},
\end{align} for $a>0$, where $\Phi(x)= \frac{1}{\sqrt{2\pi}}\int_{-\infty}^x e^{-t^2/2}\, dt$ and $\phi=\Phi'$.
\end{lemma}

The following technical lemma will be useful:

\begin{lemma} \label{lem:gauss^p}
Let $2\leq p <\infty$ and let $g_1,g_2$ be i.i.d. standard normal variables. The following properties hold:
\begin{itemize}
\item [\rm i.] The function $t\mapsto P\left ( \big| |g_1|^p -|g_2|^p \big| >t\right)$ is log-convex in $(0,\infty)$.
\item [\rm ii.] For any $r\geq 1$ we have:
\begin{align}
\left( \mathbb E \big| |g_1|^p -|g_2|^p \big|^r\right)^{1/r} \simeq r^{p/2} \sigma_p^p.
\end{align}
\end{itemize}
\end{lemma}

\noindent {\it Proof.} (i) Set $H_p(t):= P\left( \big| |g_1|^p -|g_2|^p \big| >t\right)$. Then, we may check that:
\begin{align*}
H_p(t) = \sqrt{\frac{2}{\pi}} \int_{\mathbb R} H_p(x,t) \, d\gamma_1(x),
\end{align*} where 
\begin{align*}
H_p(x,t) := \int_{(|x|^p+t)^{1/p}}^\infty e^{-y^2/2} \, dy \quad (x,t) \in \mathbb R\times (0,\infty).
\end{align*} We have the following:

\smallskip

\noindent {\it Claim 1.} For fixed $x\in \mathbb R$, the map $t\mapsto H_p(x,t)$ is log-convex on $(0,\infty)$. 

\smallskip 

\noindent To this end it suffices to check that $ H_p(x,t) \geq (H'_p(x,t))^2/ H_p''(x,t)$ for all $t>0$, equivalently:
\begin{align*}
\int_{(|x|^p +t)^{1/p} }^\infty e^{-y^2/2}\, dy \geq  \exp \left(-\frac{1}{2}(|x|^p+t)^{2/p} \right) \frac{(|x|^p+t)^{1/p} }{p-1+ (|x|^p+t)^{2/p}}.
\end{align*} The latter follows by \eqref{eq:gordon-1} (for $a=(|x|^p+t)^{1/p}$) in Lemma \ref{lem:Gordon-ineq} .

The first assertion now follows by H\"{o}lder's inequality.

\smallskip

\noindent (ii) The upper estimate is a consequence of the triangle inequality and the fact that $\sigma_{pr}^p \simeq r^{p/2} \sigma_p^p$ 
(see estimate \eqref{eq: 2.4}). 
For the lower bound we have to elaborate more. Using polar coordinates we may write:
\begin{align*}
\mathbb E\big| |g_1|^p -|g_2|^p \big|^r = 
\frac{2^{\frac{pr}{2}+2}}{\pi} \Gamma\left(\frac{pr}{2}+1 \right) \int_0^{\pi/4} \left( \cos^p\theta -\sin^p\theta \right)^r \, d\theta.
\end{align*} We have the following:

\smallskip

\noindent {\it Claim 2.} For $r\geq 1$ we have:
\begin{align*}
\int_0^{\pi/4} \left( \cos^p\theta -\sin^p\theta \right)^r \, d\theta \gtrsim (2/3)^r /\sqrt{pr}.
\end{align*}

\smallskip

\noindent Indeed; we may write:
\begin{align*}
\int_0^{\pi/4} \left( \cos^p\theta -\sin^p\theta \right)^r \, d\theta \geq \int_0^{\pi/6} \left( \cos^p\theta -\sin^p\theta \right)^r \, d\theta
\geq \left( 1-3^{-p/2} \right)^r \int_0^{\pi/6} (\cos\theta)^{pr} \, d\theta,
\end{align*} where we have used the fact that $\sin \theta \leq 3^{-1/2} \cos \theta$ for any $\theta \in [0,\pi/6]$. Next, we have:
\begin{align*}
\int_0^{\pi/6} (\cos\theta)^{pr}  \, d\theta = \frac{1}{2}B\left( \frac{pr +1}{2}, \frac{1}{2}\right)- \int_{\pi/6}^{\pi/2} (\cos \theta)^{pr} \, d\theta 
\geq \frac{1}{2}B\left( \frac{pr+1}{2}, \frac{1}{2}\right)- \frac{2\cos^{pr+1}(\pi/6)}{pr+1}.\end{align*}  A standard approximation for the Beta
function provides:
\begin{align*}
B \left(\frac{pr+1}{2}, \frac{1}{2} \right) \simeq  (pr)^{-1/2},
\end{align*} and thus, the Claim 2 follows. 

Finally, Stirling's approximation formula yields $2^{pr/2}\Gamma(\frac{pr}{2}+1)\simeq (pr)^{1/2} (pr/e)^{pr/2}$ 
and the result follows. \prend

\subsection{Functional inequalities on Gauss' space} 

First we refer to the logarithmic Sobolev inequality. In general, if $\mu$ is 
a Borel measure on $\mathbb R^n$ it is said that $\mu$ satisfies a {\it log-Sobolev inequality with constant $\rho$} if for any smooth 
function $f$ we have:
\begin{align*}
 {\rm Ent}_\mu (f^2):= \mathbb E_\mu(f^2\log f^2) - \mathbb E_\mu f^2 \log(\mathbb E_\mu f^2) \leq \frac{2}{\rho} \int \|\nabla f\|^2_2\, d\mu.
\end{align*} It is well known (see \cite{Led}) that the standard $n$-dimensional Gaussian measure $\gamma_n$ satisfies the log-Sobolev inequality 
with $\rho =1$. The next lemma, based on the classical Herbst's argument, is a useful estimate which holds for any measure 
satisfying a log-Sobolev inequality:

\begin{lemma} \label{lem:log-sob-moms}
Let $\mu$ be a measure satisfying the log-Sobolev inequality with constant $\rho>0$. 
Then, for any Lipschitz\footnote{Recall that for a Lipschitz map $f :(X,\rho)\to \mathbb R$ on some metric space
$(X,\rho)$ the Lipschitz constant of $f$ is defined by $\|f\|_{\rm Lip}=\sup_{x,y\in X, x\neq y}\frac{|f(x)-f(y)|}{\rho(x,y)}$.} 
map $f$ and for any $2\leq p< q$ we have:
\begin{align} \label{eq: 2.19}
\|f\|_{L_q(\mu)}^2- \|f\|_{L_p(\mu)}^2 \leq \frac{\|f\|_{\rm Lip}^2}{\rho}(q-p).
\end{align} In particular, we have:
\begin{align}\label{eq:2.7}
\frac{ \|f\|_{L_q(\mu)} }{\|f\|_{L_2(\mu)} }\leq \sqrt{1+\frac{q-2}{ \rho k(f)} },
\end{align} for $q\geq 2$ where $k(f):= \|f\|_{L_2(\mu)}^2/ \|f\|_{\rm Lip}^2$. Furthermore, 
\begin{align}
\frac{\|f\|_{L_2(\mu)} }{\|f\|_{L_p(\mu)} } \leq \exp\left( \frac{1/p-1/2}{ \rho k(f)}\right),
\end{align} for $0<p \leq 2$.
\end{lemma} 

\noindent {\it Proof.} The proof of the first estimate is essentially contained in \cite{SV}. The second one 
is direct application of the first for $p=2$. For the last assertion, note that by Lyapunov's convexity theorem (see \cite{HLP}) the map
$p \stackrel{\phi} \mapsto \log\|f\|_p^p$ is convex. Moreover, we have:
$p\phi'(p)-\phi(p) =\frac{{\rm Ent}_\mu(|f|^p)}{ \int |f|^p \, d\mu}$. Hence, for any $0<p<2$, the convexity of $\phi$ and 
the log-Sobolev inequality yield: 
\begin{align*}
2 \frac{\phi(2)-\phi(p)}{2-p} \leq 2\phi'(2) = \frac{{\rm Ent}_\mu(f^2)}{\|f\|_2^2} +\phi(2) \leq \frac{2}{2\rho k} +\phi(2),
\end{align*} where $k\equiv k(f)$. The result follows. \prend

\smallskip

\begin{note} When $f$ is a Lipschitz map with $k(f)\gtrsim 1$, the above two estimates imply
\begin{align} \label{eq:2.11}
\frac{\|f\|_{L_q(\gamma_n)} }{ \|f\|_{L_1(\gamma_n)} } \leq \sqrt{ 1+c_1 \frac{q-1}{k(f)}}, \quad q\geq 1.
\end{align} In the case $A$ is a centrally symmetric convex body on $\mathbb R^n$, integration in polar coordinates yields:
\begin{align}
I_r(\gamma_n,A) = c_{n,r} M_r(A),
\end{align} where 
$c_{n,r}: = \sqrt{2}[\Gamma(\frac{n+r}{2})/ \Gamma (\frac{n}{2})]^{1/r}$ and $M^r_r(A) :=\int_{S^{n-1}} \|\theta\|_A^r \, d\sigma(\theta)$.
Applying this for $A=B_2^n$ we readily see that $c_{n,r}=I_r(\gamma_n,B_2^n)$. Therefore, for $-n< s < r$ we obtain:
\begin{align} \label{eq:2.15}
\max\left\{ \frac{M_r(A)}{M_s(A)}, \frac{I_r(\gamma_n, B_2^n)}{ I_s(\gamma_n,B_2^n)}  \right\}
\leq \frac{M_r(A) I_r(\gamma_n, B_2^n)}{M_s(A) I_s(\gamma_n,B_2^n)} = \frac{I_r(\gamma_n,A)}{I_s(\gamma_n,A)} . 
\end{align}  It follows that:
\begin{align} \label{eq:iso-M}
M_q(A)/M_1(A) \leq \sqrt{ 1+c_1\frac{q-1}{k(A)} }, \quad q\geq 1.
\end{align} This estimate improves considerably upon the estimate presented 
in \cite[Statement 3.1]{LMS} or \cite[Proposition 1.10, (1.19)]{Led} in the range $1\leq q \leq k(A)$. 
For a purely probabilistic approach of this fact we refer the reader
to \cite{PPV}. 

\end{note}

It is immediate that 
\begin{align*}
\| f\|_{L_r(\gamma_n)} \lesssim \left\{ 
\begin{array}{ll} 
\displaystyle \|f\|_{L_1(\gamma_n)} , & 1\leq r\leq k(f) \\
\displaystyle \sqrt{ \frac{r}{k(f)} } \|f\|_{L_1(\gamma_n)}, & r\geq k(f)
\end{array} \right.
\end{align*} for any Lipschitz function $f$ in $(\mathbb R^n,\gamma_n)$. In \cite{LMS} it is proved that for norms 
this estimate can be reversed:

\begin{lemma} \label{lem:equiv-r-means}
Let $\|\cdot\|_A$ be a norm on $\mathbb R^n$. Then, we have:
\begin{align*}
I_r(\gamma_n,A) \simeq \left\{ \begin{array}{ll}
I_1(\gamma_n,A) , & r\leq k(A) \\
\displaystyle \sqrt{\frac{r}{k(A)} } I_1(\gamma_n, A), & r\geq k(A)
\end{array} \right. .
\end{align*}
\end{lemma}

\medskip

This result implies the next well known fact:

\begin{proposition} \label{prop: sharp-tails}
Let $\|\cdot \|\equiv \|\cdot\|_A$ be a norm on $\mathbb R^n$. Then, we have:
\begin{align*}
c\exp(-Ct^2k) \leq P \left(  \|X\| >(1+t) \mathbb E\|X\| \right) \leq C\exp(-ct^2 k),
\end{align*} for $t\geq 1$. Moreover, one has:
\begin{align*}
\left( \mathbb E \big| \|X\| -\mathbb E\|X\|\big|^r \right)^{1/r} \simeq \sqrt{\frac{r}{k}} \mathbb E\|X\|,
\end{align*} for all $r\geq k$, where $k\equiv k(A)$ and $X$ is a standard Gaussian $n$-dimensional random vector.
\end{proposition}

\noindent {\it Sketch of proof of Proposition \ref{prop: sharp-tails}.} Set $I_r\equiv I_r(\gamma_n,A)$. There exists $c_1 \in (0,1)$ such that $I_s\geq c_1\sqrt{s/k}I_1$ for all $s>k$ by Lemma \ref{lem:equiv-r-means}. Thus, for $t\geq 1$, 
if we choose $r>k$ by 
$c_1\sqrt{r/k}=4t$, we may write:
\begin{align*}
P \left( \|X\| >\frac{1}{2} I_r \right) \leq P \left( \|X\| > \frac{c_1}{2} \sqrt{r/k} I_1 \right) \leq P(\|X\| \geq (1+t)I_1).
\end{align*} On the other hand the Paley-Zygmund inequality (Lemma \ref{lem: PZ-ineq}) yields:
\begin{align*}
P \left( \|X\| >\frac{1}{2} I_r \right) \geq (1-2^{-r})^2 (I_r/I_{2r})^{2r} \geq c_2e^{-C_2r} \geq c_2\exp(-C_2' t^2 k), 
\end{align*} where we have also used the fact that $I_r\simeq I_{2r}$ which follows by Lemma \ref{lem:equiv-r-means}. 
For the second assertion we apply integration by parts and we use the first estimate. \prend

\medskip

The above estimate shows that the large deviation estimate for norms with respect to $\gamma_n$ is completely settled. 
Therefore for the concentration inequalities we are interested in, we may restrict ourselves to the range $0<\varepsilon<1$.

Other important functional inequalities related to the concentration of measure phenomenon are the Poincar\'{e} inequalities. 
Using a standard variational argument (see \cite{Led}) one can show that any measure which satisfies a log-Sobolev 
inequality with constant $\rho$ also satisfies a Poincar\'{e} inequality with constant $\rho$, i.e.
\begin{align} \label{eq:Poin-ineq}
 \rho {\rm Var}_\mu(f) \leq \int_{\mathbb R^n} \|\nabla f\|_2^2 \, d\mu,
\end{align} for any smooth function $f$.

A refinement of the Poincar\'{e} inequality was proved by Talagrand in \cite{Tal} for the 
discrete cube $\{-1,1\}^n$ (see also \cite{BLM} for a recent exposition) and its continuous version, in the Gaussian context,
was presented in \cite{CL} (see also \cite{Cha}): 

\begin{theorem}[Talagrand's $L_1-L_2$ bound] \label{thm:Talagrand bd} Let $f: \mathbb R^n \to \mathbb R$ be a smooth function. If $A_i:=\|\partial_i f\|_{L_2(\gamma_n)} $ and $a_i:= \|\partial_i f\|_{L_1(\gamma_n)}$, then one has:
 \begin{align*}
 {\rm Var}_{\gamma_n} (f) \leq C \sum_{i=1}^n \frac{A_i^2}{1+\log ( A_i / a_i) },
 \end{align*} where $\partial_i f$ stands for the partial derivative $\partial f / \partial x_i$.
 \end{theorem}
 
\noindent This inequality will be used in order to prove concentration for the $\ell_p$ norm when $p$ is sufficiently large.

Pisier discovered in \cite{Pis} another Gaussian inequality which contains the $(r,r)$-Poincar\'{e} inequalities and the 
Gaussian concentration inequality as a special case (see Remarks \ref{rems:2-12}). 

\begin{theorem}\label{thm:Pis-ineq}
Let $\phi:\mathbb R\to \mathbb R$ be a convex function and let $f:\mathbb R^n\to \mathbb R$ be $C^1$-smooth. Then, if $X,Y$ are
independent copies of a Gaussian random vector, we have:
\begin{align*} 
\mathbb E \phi\left( f(X)-f(Y)\right)\leq  \mathbb E \phi \left( \frac{\pi}{2} \langle \nabla f(X), Y\rangle \right).
\end{align*}
\end{theorem}

\begin{remarks} \label{rems:2-12} \rm
\noindent 1. {\it $(r,r)$-Poincar\'{e} inequalities.} For $\phi(t)=|t|^r, \, r\geq 1$ we get:
\begin{align} \label{eq: 2.9}
\|f-\mathbb Ef\|_{L_r(\gamma_n)} \simeq \left(\mathbb E |f(X)-f(Y)|^r \right)^{1/r} \leq \frac{\pi}{2}\sigma_r \left( \mathbb E \|\nabla f(X)\|_2^r \right)^{1/r}.
\end{align} In particular for $r=2$ we have ${\rm Var}(f(X)) \leq \frac{\pi^2}{8} \mathbb E \|\nabla f(X)\|_2^2$, which is the Gaussian 
Poincar\'e inequality with non-optimal constant.

\noindent 2. {\it Gaussian concentration.} The choice $\phi_\lambda(t)=\exp(\lambda t), \; \lambda>0 $ and a standard optimization argument on $\lambda$ (see \cite{Pis} for the details) yield:
\begin{align} \label{eq:2.14}
P( |f(X)- \mathbb Ef(X) |>t )\leq 2 \exp(-t^2/(2\pi^2 \|f\|^2_{\rm Lip})),
\end{align} for all $t>0$. Alternatively, we may conclude a similar estimate by equations \eqref{eq: 2.9} and Markov's inequality.
\end{remarks}

\subsection{Negative moments of norms} 

The next result is due to Klartag and Vershynin from \cite{KV} (see also
\cite{LO} for a similar estimate as \eqref{eq:gauss-sb} with $k(A)$ instead of $d(A)$):

\begin{proposition} \label{prop:small-ball}
Let $A$ be a centrally symmetric convex body on $\mathbb R^n$. We define:
\begin{align}
d(A):= \min\left\{ n, -\log \gamma_n \left( \frac{m}{2} A \right)\right\},
\end{align} where $m$ is the median of $x \mapsto \|x\|_A$ with respect to $\gamma_n$. Then, one has:
\begin{align} \label{eq:gauss-sb}
\gamma_n\left( \left \{ x: \|x\|_A \leq c\varepsilon \mathbb E\|X\|_A \right \} \right)\leq (C\varepsilon)^{cd(A)},
\end{align} for all $0<\varepsilon<\varepsilon_0$ where $\varepsilon_0>0$ is an absolute constant. 
Moreover, for all $0<k< d(A)$ we have:
$ I_{-k}(\gamma_n, A) \geq c I_1(\gamma_n, A)$. Note that $d(A)> c_1 k(A)$.
\end{proposition}

Note that this result implies that the negative moments exhibit stable behavior up to the point $d(A)$. 
However, one can show that up to the critical dimension the moments of any norm with respect to the
Gaussian (or the uniform on the sphere) measure are almost constant, thus complementing the estimates 
\eqref{eq:2.11} and \eqref{eq:iso-M}. In order to quantify the latter we need the next consequence of 
Proposition \ref{prop:small-ball}.

\begin{lemma}
\label{lem:reduc-neg-moms} Let $A$ be a centrally symmetric convex 
body on $\mathbb R^n$ which satisfies the small ball probability estimate:
\begin{align*}
\gamma_n( \varepsilon I_1 A) < (K\varepsilon)^{\alpha d},
\end{align*} for all $0<\varepsilon< \varepsilon_0$ ($K, \alpha>0$). Then, for all $r,s>0$ with $r+s< \alpha d/3$ we have:
\begin{align*}
I_{-r-s}^{-r-s}(\gamma_n, A) \leq \left(\frac{CK}{I_1}\right)^s I_{-r}^{-r}(\gamma_n, A),
\end{align*} where $I_1=I_1(\gamma_n,A)$ and $C>0$ is an absolute constant. 
\end{lemma}

\noindent {\it Proof.} We set $I_q=I_q(\gamma_n,A)$. For any $0<\varepsilon <\varepsilon_0$ we may write:
\begin{align*}
I_{-r-s}^{-r-s} = \int \frac{1}{\|x\|_A^{r+s}} \, d\gamma_n(x) &\leq
\frac{1}{(\varepsilon I_1)^s} \int \frac{1}{\|x\|_A^r} \, d\gamma_n(x)+\int_{\varepsilon I_1A} \frac{1}{\|x\|_A^{r+s}} \, d\gamma_n(x) \\
&\leq \frac{1}{(\varepsilon I_1)^s} I_{-r}^{-r} + (K\varepsilon)^{\alpha d/2} I_{-2(r+s)}^{-r-s},
\end{align*} by the Cauchy-Schwarz inequality. Note that the small ball probability assumption implies that: $I_{-s}\geq c\varepsilon_0 I_1$
for all $0<s< 2\alpha d/3$. Thus, if $r+s<\alpha d/3$ we get $I_{-2(r+s)} >c_1 I_{-(r+s)}$ and previous estimate yields:
\begin{align*}
I_{-r-s}^{-r-s} <\frac{1}{(\varepsilon I_1)^s} I_{-r}^{-r} + (K\varepsilon)^{\alpha d/2}  c_1^{-r-s} I_{-r-s}^{-r-s}.
\end{align*} Choosing $\varepsilon$ small enough so that $(K\varepsilon)^{\alpha d/2}< c_1^{r+s}/2$, say $0< \varepsilon \leq c_1/(2K)$, 
we conclude the result. \prend

\begin{theorem}\label{thm:stability-moms}
Let $A$ be a centrally symmetric convex body on $\mathbb R^n$. Then, one has:
\begin{align*}
\frac{I_r(\gamma_n, A)}{I_{-r}(\gamma_n,A)} \leq 1+\frac{Cr}{k(A)},
\end{align*} for all $0<r<ck(A)$, where $C,c>0$ are absolute constants.
\end{theorem}

\noindent {\it Proof.} We present the argument in two steps:

\smallskip

\noindent {\it Step 1.}  (positive moments). We use the log-Sobolev inequality to estimate the growth of the moments. 
The basic observation is that:
\begin{align*}
\frac{d}{dr} \left( \log \|f\|_{L_r(\mu)} \right) = \frac{{\rm Ent}_\mu (|f|^r)}{r^2 \|f\|_{L_r(\mu)}^r},
\end{align*} for any Lipschitz function $f$. Apply this for the function $f=\|\cdot\|_A $ to get:
\begin{align*}
(\log I_r)' \leq \frac{1}{2I_r^r} \mathbb E \|X\|_A^{r-2} \| \, \nabla \|X\|_A \, \|_2 \leq \frac{b^2}{2I_r^r} I_{r-2}^{r-2},
\end{align*} for all $r>0$, where $b=b(A)$ the Lipschitz constant of $\|\cdot\|_A$. It is easy to see that $(\log I_r)' \leq \frac{1}{2k(A)}$ for $r\geq 2$, 
while for $0<r<2$ we may write:
\begin{align}\label{eq:4.11}
(\log I_r)' \leq \frac{b^2}{2 I_{-(2-r)}^2}\leq \frac{C_1b^2}{I_1^2} \leq \frac{C_1'}{k(A)},
\end{align} where we have used Proposition \ref{prop:small-ball}. Using \eqref{eq:4.11} we may write:
\begin{align} \label{eq:2.41}
\log (I_r/I_0)= \int_0^r (\log I_t)'\, dt \leq \int_0^r \frac{C_1}{k} \, dt=\frac{C_1r}{k},
\end{align} for all $r>0$.
\smallskip

\noindent {\it Step 2.} (negative moments). As before, using the log-Sobolev inequality, for all $0<r<c_1d(A)$ we may write:
\begin{align*}
(\log I_{-r})' \geq -\frac{b^2}{2I_{-r}^{-r}}I_{-r-2}^{-r-2} \geq -\frac{C_2 b^2}{I_1^2} \geq -\frac{C_2'}{k(A)},
\end{align*} where we have used Lemma \ref{lem:reduc-neg-moms}. The same reasoning applied to \eqref{eq:2.41} shows that
$\log (I_{-r}/I_0)\geq -\frac{C_2r}{k}$, for all $0<r<c_1d(A)$. Combining the two steps and restricting to $0<r<c_2k(A)$ we conclude 
the result. \prend


\section{The Gaussian variance of the $\ell_p$ norm}

A standard method for bounding the variance is the concentration 
inequality \eqref{eq:2.14}, e.g. see \cite{LMS} or \cite[Proposition 1.9]{Led}. An integration by 
parts argument implies that if $f:\mathbb R^n\to \mathbb R$ is $L$-Lipschitz function, then
${\rm Var}(f) \lesssim L^2$. In particular, if $f(x)=\|x\|_p$ this estimate yields:
\begin{align*}
{\rm Var} \|X\|_p \lesssim b^2(B_p^n) \simeq \max\{ n^{2/p-1}, 1\}, \quad 1\leq p\leq \infty.
\end{align*}
For $1\leq p\leq 2$ this estimate turns out to be the correct one. But, for $2<p\leq \infty$ this method gives bounds
which are far from the actual ones. The purpose of this Section is to compute the correct order of magnitude 
for the Gaussian variance of the $\ell_p$ norm. Our first approach lies in determining the limit distribution of 
the sequence of variables $(\|g\|_{\ell_p^n})_{n=1}^\infty$. Here $\|g\|_{\ell_p^n}$ stands for the $\ell_p$ norm of the 
$n$-dimensional ``truncation" of the sequence $(g_i)_{i=1}^\infty$ of i.i.d. standard Gaussian random variables, 
i.e. $\|g\|_{\ell_p^n}:=(\sum_{i\leq n}|g_i|^p)^{1/p}$.

\subsection{ The variance of the $\ell_p$ norm for $1 \leq p <\infty$}  

In this case we use the next Proposition known in Statistics as the "Delta Method" (for a proof see \cite{Cra}):

\begin{proposition} \label{prop: delta-meth} Let $\theta, \sigma \in \mathbb R$ and let $(Y_n)$ be a sequence of random 
variables that satisfies $n^{1/2}(Y_n-\theta) \longrightarrow N(0,\sigma^2)$ in
distribution. For the differentiable function $h$ assume that $h'(\theta)\neq 0$. Then,
\begin{align*}
n^{1/2} (h(Y_n)-h(\theta)) \longrightarrow N(0,\sigma^2(h'(\theta))^2)
\end{align*} in distribution. 
\end{proposition}

Now we may prove the next asymptotic estimate:

\begin{theorem} \label{thm:up-low bound var}
Let $1\leq p <\infty$. Let $(\xi_j)_{j=1}^\infty$ be sequence of i.i.d random variables with $m_{3p}^{3p}:=\mathbb E|\xi_1 |^{3p}<\infty$.
Then, there exist positive constants $c_p,C_p$ depending only on $p$ and the distribution of $(\xi_j)$ such that:
\begin{align*}
c_p n^{\frac{2}{p}-1} \leq {\rm Var}\| \xi \|_{\ell_p^n} \leq C_p n^{\frac{2}{p}-1},
\end{align*} for all $n$, where $\|\xi\|_{\ell_p^n}^p =\sum_{j\leq n} |\xi_j|^p$. 
\end{theorem}

\noindent {\it Proof.} Let $Y_n:=\frac{1}{n}\sum_{j=1}^n |\xi_j|^p$. Then by the Central Limit Theorem 
we know that:
\begin{align*}
\sqrt{n}(Y_n-m_p^p) \longrightarrow N(0, v_p^2)
\end{align*} in distribution, where $v_p^2 := {\rm Var}|\xi_1|^p$. Consider the function $h(t)=t^{1/p}, \; t>0$ and
apply Proposition \ref{prop: delta-meth} to get:
\begin{align*}
\zeta_n:= \sqrt{n}(n^{-1/p}\|\xi\|_p-m_p) \longrightarrow N \left( 0, \frac{v_p^2}{p^2} m_p^{2(1-p)} \right),
\end{align*} in distribution. Using the fact that $m_{3p}<\infty$ we may conclude the uniform integrability of $(\zeta_n^2)_{n=1}^\infty$:

\noindent {\it Claim.} For all $n\geq 1$ we have:
\begin{align*}
\mathbb E |\zeta_n|^3 \lesssim m_{3p}^{3p}/ m_p^{3(p-1)}.
\end{align*}

\noindent {\it Proof of Claim.} We may write:
\begin{align*}
\mathbb E |\zeta_n|^3 & =n^{\frac{3}{2}-\frac{3}{p}} \mathbb E \left| \|\xi\|_p -n^{1/p} m_p \right|^3
\leq \frac{ n^{\frac{3}{2}-\frac{3}{p}}  }{ (n^{1/p}m_p)^{3(p-1)} }  \mathbb E \left| \|\xi\|_p^p -nm_p^p \right|^3 
\leq \frac{ n^{-3/2}  }{ m_p^{3(p-1)} } \mathbb E \left | \|\xi\|_p^p - \|\xi'\|_p^p \right|^3,
\end{align*} where $\xi'$ is an independent copy of $\xi$ and we have also used the numerical inequality $a^{p-1} |z-a| \leq |z^p-a^p|$ for $z\geq 0, a>0, p\geq 1$ and Jensen's inequality. Finally, a standard symmetrization argument yields:
\begin{align*}
\mathbb E \left| \sum_{j=1}^n \left( |\xi_j|^p -|\xi_j'|^p \right) \right|^3 \lesssim \mathbb E \left[\sum_{j=1}^n \left( |\xi_j|^p -|\xi_j'|^p \right)^2 \right]^{3/2} \lesssim n^{3/2} \mathbb E \left| |\xi_1|^p -|\xi_1'|^p \right|^3 \lesssim n^{3/2} m_{3p}^{3p},
\end{align*} where  we have also used Jensen's inequality, again. This proves the claim. 

Hence, we may conclude:
\begin{align} \label{eq:3.7}
n^{1-\frac{2}{p}} {\rm Var}(\| \xi \|_p) =
{\rm Var} \left( n^{\frac{1}{2}-\frac{1}{p}} \|\xi \|_p \right) = 
{\rm Var} \left[ \sqrt{n}(n^{-1/p}\| \xi \|_p-m_p) \right] \to \frac{v_p^2}{p^2}m_p^{2(1-p)},
\end{align}  as $n\to \infty$ and the result follows. \prend

\smallskip

\begin{remark} \rm The reader should notice that, for fixed $p \geq 1$, the dependence we obtain on the dimension is the same regardless the randomness we choose for the underlying variables $(\xi_i)$. In addition the argument is essentially based on the 
stochastic independence. 
Moreover, in the case that $(\xi_i)$ are standard normals, the above limit value is estimated as: 
\begin{align} \label{eq:3.8}
\frac{v_p^2}{p^2}m_p^{2(1-p)} \sim  \frac{1}{e\sqrt{2}} \frac{2^p}{p},  \quad 
p\to \infty.
\end{align} This suggests that the constants $c_p,C_p$ should depend exponentially on $p$. 
\end{remark}

\subsection {The variance of the $\ell_\infty$ norm} 

Of course the variance in that case can be computed by employing the tail estimates
for the $\ell_\infty$-norm proved in \cite{Sch2}. However, we prefer here to give a proof of a more 
"probabilistic flavor". Actually, the argument we present below works for all i.i.d. random variables with exponential tails, 
but we shall focus on Gaussians.
Let $(g_i)_{i=1}^\infty$ be independent, standard Gaussian random variables and let $Y_n:= \max_{i\leq n} |g_i|, \; n\geq 2$. 
We set $a_n:= -\Phi^{-1}(\frac{1}{2n}) >0$. Note that $a_n\to \infty$ and Gordon's inequality \eqref{eq:gordon-2} 
shows that $a_n \sim \sqrt{2\log n}$ as $n\to \infty$. We define $W_n: = a_n(Y_n-a_n)$ and we have the next well known fact 
(see \cite[\S \, 9.3]{Dav}):

\begin{proposition}
Let $\eta$ be a Gumbel random variable, that is the cumulative distribution function of $\eta$ is given as: 
\begin{align*}
F_\eta(t) := \exp(-e^{-t}), \; t\in \mathbb R. 
\end{align*} If $(W_n)$ is the sequence defined above, then for every $t\in \mathbb R$ we have:
\begin{align*} 
\mathbb P (W_n\leq t) \to \exp(-e^{-t}) ,
\end{align*} that is $W_n$ converges to the Gumbel variable in distribution.
\end{proposition}

For the random variable $\eta$ it is known that $\mathbb E (\eta) =\gamma$ (the Euler-Mascheroni constant) and 
${\rm Var} (\eta) = \pi^2/6$. Therefore, we obtain:
\begin{align*}
a_n^2 {\rm Var}(Y_n) = {\rm Var}(W_n) \to {\rm Var} (\eta) ,
\end{align*} as $n\to \infty$. This proves the following:

\begin{theorem}\label{thm:var-ell-infty} If $Z$ is an $n$-dimensional standard Gaussian random vector, we have:
\begin{align*}
{\rm Var} \|Z\|_\infty ={\rm Var}_{\gamma_n} \|x\|_\infty \simeq (\log n)^{-1}.
\end{align*} 
\end{theorem}
 
\medskip

It should be noticed that the dependence on dimension we get for fixed $1\leq p<\infty$ is polynomial in $n$ while
for $p=\infty$ is logarithmic in $n$. 
As we have already explained this ``skew'' behavior relies on the fact that as $p$ grows, the constants in the equivalence 
should be expected to be exponential in $p$ (see \eqref{eq:3.7} and \eqref{eq:3.8}). In the rest of the paragraph we try to 
study and quantify this phenomenon. Our aim is to give as sharp bounds as possible and describe the behavior of $p$ along $n$, too.

\subsection{Tightening the bounds} 

The purpose of this subsection is to provide continuous bounds in terms of 
$p$ for the variance of the $\ell_p$ norm when dimension $n \to \infty$ and $p$ varies from $1$ to $\infty$ (along with $n$).
One can easily see that:
\begin{align*}
c_1 p \leq n^{1-2/p} {\rm Var}\|X\|_p \leq c_2 p {\rm Var}|g_1|^p \simeq p(2p/e)^p,
\end{align*}
by just comparing with the variance of the $\ell_2$ norm and the $p$-th power of the $\ell_p$ norm. 
Below, we show that one can always have better estimates. In order to prove these estimates we will use the following:

\begin{lemma} \label{lem:bd-var} Let $4\leq p\leq \infty$. Then one has:
\begin{align*}
I_r (\gamma_n,B_p^n) /  I_{-r} (\gamma_n,B_p^n) \leq \exp \left( \frac{C_1r}{k_{p,n} \log n} \right) , \quad 0<r< c_1\sqrt{k_{p,n} \log n},
\end{align*}
where $k_{p,n}\equiv k(B_p^n)$.
\end{lemma} 

We postpone the proof of this Lemma to Section 4 (Theorem \ref{thm:stability-r-means}).

\subsubsection{Upper bound (via Talagrand's inequality)} 

For $p>1$ we have: 
$\partial_i \|x\|_p =\frac{|x_i|^{p-1}}{\|x\|_p^{p-1}} \sgn{x_i}$ a.s. Thus, one has:
\begin{align*}
A^2:= \big \| \partial_i \| \cdot \|_p \big\|_{L_2}^2 \leq \sigma_{2p-2}^{2p-2}I_{-2(p-1)}^{-2(p-1)}(\gamma_{n-1},B_p^{n-1}), \; 
a: = \big \| \partial_i \|\cdot \|_p \big\|_{L_1} \leq \sigma_{p-1}^{p-1}I_{-(p-1)}^{-(p-1)}(\gamma_{n-1},B_p^{n-1}).
\end{align*} Set $I_s(\gamma_{n-1},B_p^{n-1})\equiv I_s$. Thus, direct application of Theorem \ref{thm:Talagrand bd} yields:
\begin{align} \label{eq:6}
{\rm Var}(\|X\|_p)\leq C n\ \frac{\sigma_{2p-2}^{2p-2} I_{-2(p-1)}^{-2(p-1)}  }{1+\log\left(\frac{\sigma_{2p-2}^{p-1} }{\sigma_{p-1}^{p-1}} \frac{I_{-(p-1)}^{p-1}}{I_{-2(p-1)}^{p-1} } \right) } \leq C_1 n\frac{\sigma_{2p-2}^{2p-2}/I_{-2(p-1)}^{2(p-1)} }{p},
\end{align} where we have used the fact that $\left(\sigma_{2p-2}/ \sigma_{p-1}\right)^{p-1} \simeq 2^p$, which 
follows by \eqref{eq: 2.4}. As long
as $2p<c_1\sqrt{k_{p,n} \log n}$, which is satisfied when $p\leq c_0\log n$ for some sufficiently small
absolute constant $c_0>0$ in view of Proposition \ref{prop:mean-ell-p}, we may apply Lemma \ref{lem:bd-var} to get:
\begin{align*}
I_{-2(p-1)}^{2(p-1)} \geq e^{ - \frac{c' p^2}{k_{p,n} \log n}}I_p^{2(p-1)} \geq c_1' \sigma_p^{2(p-1)} (n-1)^{2-2/p}. 
\end{align*} Plug this estimate in \eqref{eq:6} we derive the upper bound:
\begin{align*}
{\rm Var}\|X\|_p \leq C_2 \frac{ \sigma_{2p-2}^{2p-2} }{\sigma_p^{2p-2} p} n^{\frac{2}{p}-1} \simeq \frac{2^p}{p} n^{2/p-1}.
\end{align*} Note that this is exactly of the same order as the one we obtained at the limit value using the Delta Method.

\subsubsection{Lower bound (via Talagrand's inequality)}  

Here we will use the next numerical result:

\begin{lemma} \label{lem:2-sided-ineq}
Let $a,b>0$ and $0<\theta \leq 1$. Then, we have:
\begin{align*}
\theta |a-b| \left( \frac{2}{a+b}\right)^{1-\theta} \leq |a^\theta -b^\theta| \leq \theta |a-b| \frac{a^{\theta-1} + b^{\theta-1}}{2} .
\end{align*} \end{lemma}

\noindent {\it Proof.} We may assume without loss of generality that $0<a<b$ and $0<\theta<1$. If we set $f(t)=t^{\theta -1}, \; t>0$, 
note that $f$ is convex in $[a,b]$, hence the estimate follows by the Hermite-Hadamard inequality (see \cite{HLP}). \prend

\medskip

Applying the lower bound of Lemma \ref{lem:2-sided-ineq} for $a=\|X\|_p^p, \; b=\|Y\|_p^p$ and $\theta =1/p$, where 
$X,Y$ are independent and $X,Y \sim N({\bf 0}, I_n)$, we obtain:
\begin{align} \label{eq:10}
2 {\rm Var}\|X\|_p = \mathbb E( \|X\|_p-\|Y\|_p)^2 \geq \frac{2^{2/q}}{p^2} \mathbb E \frac{(\|X\|_p^p -\|Y\|_p^p )^2}{( \|X\|_p^p +\|Y\|_p^p)^{2/q}}
\geq \frac{1}{p^2} \mathbb E \left|\sum_{i=1}^n \frac{|X_i|^p-|Y_i|^p}{S^{1/q}}\right|^2,
\end{align} where $q$ is the conjugate exponent of $p$, i.e. $1/p+1/q=1$ and 
\begin{align*}
S:= \|X\|_p^p+\|Y\|_p^p = \|Z\|_p^p , \quad Z=(Z_1, \ldots ,Z_{2n}) \sim N({\bf 0},I_{2n}).
\end{align*} 
Now we observe that the variables $\eta_j : =\frac{|X_j|^p - |Y_j|^p}{S^{1/q}}$ have the same distribution and 
satisfy $\mathbb E(\eta_i \eta_j)=0$ for $i\neq j$.
Therefore, we have:
\begin{align*}
\mathbb E \left|\sum_{i=1}^n \frac{|X_i|^p-|Y_i|^p}{S^{1/q}}\right|^2 = 
\mathbb E \left|\sum_{i=1}^n \eta_i\right|^2 =\sum_{i=1}^n \mathbb E \eta_i^2 = n\mathbb E \eta_1^2
\end{align*}
Hence, estimate \eqref{eq:10} becomes:
\begin{align*}
{\rm Var}\|X\|_p \geq 
\frac{n}{2p^2} \mathbb E \frac{(|X_1|^p-|Y_1|^p)^2}{S^{2/q}} = 
\frac{n}{p^2} \left( \mathbb E \frac{|X_1|^{2p}}{S^{2/q}} - \mathbb E \frac{|X_1|^p |Y_1|^p}{S^{2/q}}\right).
\end{align*} Let  $T:=\sum_{i>1}|X_i|^p+ \sum_{i>1}|Y_i|^p$. Note that $T\leq S$, thus we obtain:
\begin{align} \label{eq:var-1} 
{\rm Var} \|X\|_p \geq \frac{n}{p^2} \left[ \mathbb E\frac{|Z_1|^{2p}}{S^{2/q}}-\sigma_p^{2p} \mathbb E (T^{-2/q}) \right].
\end{align} 
An application of Lemma \ref{lem:bd-var} yields 
\begin{align} \label{eq:var-2}
\mathbb E(T^{-2/q}) \lesssim \frac{1}{\sigma_p^{2p-2}(n-1)^{2-2/p}},
\end{align} as long as $p\leq c_0\log n$. For the term $\mathbb E\frac{|Z_1|^{2p}}{S^{2/q} }$ we may write:
\begin{align*}
\mathbb E\frac{|Z_1|^{2p}}{S^{2/q} } = (2n)^{-1} \mathbb E \frac{\|Z\|_{2p}^{2p}}{S^{2/q}} 
=(2n)^{-1} \mathbb E \frac{\|Z\|_{2p}^{2p}}{\|Z\|_p^{2(p-1)}} \geq \frac{(\mathbb E \|Z\|_{2p}^p)^2}{2n \mathbb E\|Z\|_p^{2(p-1)} },
\end{align*} where we have used that the variables $|Z_j|^{2p}/S^{2/q}$ are equidistributed and the Cauchy-Schwarz inequality. 
Now by using Lemma \ref{lem:bd-var} again we obtain:
\begin{align*}
\mathbb E\|Z\|_{2p}^{2p} \leq e^{ \frac{cp^2}{k_{2p,2n} \log n}} (\mathbb E \|Z\|_{2p}^p)^2 \leq C_1 (\mathbb E \|Z\|_{2p}^p)^2
\end{align*} and similarly we have: $ \mathbb E \|Z\|_p^{2(p-1)} \leq C_2 (\mathbb E \|Z\|_p^p)^{2(p-1)/p}$, as long as $p\leq c_0\log n$. Therefore, we get:
\begin{align} \label{eq:var-3}
\mathbb E\frac{|Z_1|^{2p}}{S^{2/q} } \geq \frac{c_3}{n} \frac{ \mathbb E \|Z\|_{2p}^{2p}}{ (\mathbb E \|Z\|_p^{p})^{2(p-1)/p} } \simeq  
\frac{\sigma_{2p}^{2p}}{ n^{2-2/p} \sigma_p^{2(p-1)}}.
\end{align}
Inserting \eqref{eq:var-2} and \eqref{eq:var-3} in \eqref{eq:var-1} we get:
\begin{align*}
{\rm Var}\|X\|_p \geq \frac{c_4n}{p^2} \left[ c_5 \frac{\sigma_{2p}^{2p} }{n^{2-2/p}\sigma_p^{2(p-1)}} - c_6\frac{\sigma_p^{2p} }{\sigma_p^{2p-2} n^{2-2/p}} \right] = \frac{c_4c_5 \sigma_p^2}{p^2 n^{1-2/p} } \left[ \frac{\sigma_{2p}^{2p} }{\sigma_p^{2p}}-\frac{c_6}{c_5}\right].
\end{align*} Taking into account that $(\sigma_{2p}/\sigma_p)^{2p}\simeq 2^p$ we may conclude:
\begin{align}
{\rm Var}(\|X\|_p) \geq c_7\frac{2^p}{p} n^{2/p-1},
\end{align} provided that $p$ is greater than some large absolute constant. 

\bigskip

Finally, for much larger values of $p$, namely for $p\geq c_0\log n$, we employ Theorem \ref{thm:Talagrand bd} again. This is 
an extension of the known argument for $\ell_\infty$, which can be found in \cite{Cha}. 
As before, if $a_i:= \|\partial_i  f\|_{L_1(\gamma_n)}$  we may write:
\begin{align*} 
a_i = \int_{\mathbb R^n} \frac{|x_i|^{p-1}}{\|x\|_p^{p-1}}\, d\gamma_n(x)= \frac{1}{n} \int_{\mathbb R^n} \left( \frac{\|x\|_{p-1}}{\|x\|_p}\right)^{p-1}\, d\gamma_n(x) \leq \frac{n^{1/p}}{n}=n^{-1/q},
 \end{align*} where in the last step we have used estimate \eqref{eq: Holder - p-norms} and $q$ is the conjugate exponent of $p$. 
 Moreover, we have:
\begin{align*} 
A_i^2:=\|\partial_i f\|_{L_2(\gamma_n)}^2 = \int_{\mathbb R^n}\frac{|x_i|^{2p-2}}{\|x\|_p^{2p-2}}= \frac{1}{n} \int_{\mathbb R^n} 
\left(\frac{\|x\|_{2p-2}}{\|x\|_p}\right)^{2p-2}\, d\gamma_n(x) \leq 1/n,
\end{align*} by the estimates \eqref{eq: Holder - p-norms} again. These bounds and Theorem \ref{thm:Talagrand bd} yield:
\begin{align} \label{eq:var-4}
{\rm Var}(\|X\|_p) \leq C  \sum_{i=1}^n \frac{A_i^2}{1+ \frac{1}{q}\log n +\log A_i } \leq \frac{C}{\log n}, \quad p\geq c_0\log n
\end{align} where we have used the monotonicity of $t\mapsto \frac{t^2}{1+ \frac{1}{q}\log n +\log t}$ and that 
$q \ll 2$.

\smallskip

Finally, let us note that the variance of the $\ell_p$ norm stabilizes around $\frac{1}{\log n }$ for $p> (\log n)^2$. This is a special
case of the next reverse concentration estimate: 
\begin{proposition} \label{prop: 3.5}
Let $p> (\log n)^2$ and let $X$ be an $n$-dimensional standard Gaussian random vector. Then we have:
\begin{align*}
P\left( \big| \|X\|_p -\mathbb E \|X\|_p \big|  > \varepsilon  \mathbb E\|X\|_p \right) \geq c e^{-C\varepsilon \log n} ,
\end{align*} for all $0<\varepsilon<1$, where $C,c>0$ are absolute constant.
In particular, we have:
\begin{align*}
{\rm Var}\|X\|_p \simeq \frac{1}{\log n}.
\end{align*}
\end{proposition}

\noindent {\it Proof.} Consider $\frac{2}{\log n} <\varepsilon <1$ and write:
\begin{align*}
P( \|X\|_p > (1+\varepsilon)\mathbb E\|X\|_p) 
& \geq P(\|X\|_\infty > (1+\varepsilon)n^{1/p} \mathbb E\|X\|_\infty) \\
& \geq P(\|X\|_\infty > (1+2\varepsilon) \mathbb E\|X\|_\infty)  > c e^{-C\varepsilon \log n},
\end{align*} where we have used \eqref{eq: Holder - p-norms} and at the last step the concentration from \cite{Sch2}. Hence,
\begin{align*}
P\left( \big| \|X\|_p -\mathbb E\|X\|_p \big| >\varepsilon \mathbb E\|X\|_p \right) & \geq c' e^{-C \varepsilon \log n} ,
\end{align*} for all $0<\varepsilon <1$. For the second assertion we may write:
\begin{align*} {\rm Var}(\|X\|_p ) &= 2 (\mathbb E\|X\|_p)^2 \int_0^{\infty} t P \left( \big| \|X\|_p -\mathbb E\|X\|_p \big| > t\mathbb E\|X\|_p 
\right) \, dt \\
& \geq 2 c'(\mathbb E\|X\|_p)^2 \int_0^1 t e^{-Ct\log n} \, dt  \gtrsim \frac{(\mathbb E\|X\|_p)^2}{(\log n)^2}.
\end{align*} The result follows by Proposition \ref{prop:mean-ell-p}. \prend

\medskip

The results of this paragraph can be summarized in the next:

\begin{theorem} \label{thm:var-ell_p}
There exist absolute constants $c_0, c_1,C_1>0$ with the following property: For all $n$ large enough and for any 
$1\leq p\leq c_0\log n$ we have:
\begin{align}
c_1\frac{2^p}{p} \leq n^{1-\frac{2}{p}} {\rm Var} \|X\|_p \leq C_1 \frac{2^p}{p}.
\end{align} If $p>c_0 \log n$ then we have:
\begin{align} 
{\rm Var}\|X\|_p \leq \frac{C_1} {\log n} ,
\end{align} whereas for $p\geq (\log n)^2$ we also have:
\begin{align}
{\rm Var} \|X\|_p \geq \frac{c_1}{\log n},
\end{align} where $X\sim N({\bf 0},I_n)$.
\end{theorem}

\begin{note} While this paper was under review, Tikhomirov \cite{Tik3} improved Proposition \ref{prop: 3.5}
by extending the range to $p \geq C_0\log n$ (his proof gives $C_0=12$). 
In particular, ${\rm Var} \|X\|_p \gtrsim (\log n)^{-1}$ for $p \geq C_0 \log n$. 
We present his argument in the Appendix. This only leaves a relatively small interval $(c_0\log n, C_0\log n)$, 
for which the behavior of the variance is not exactly determined. In other words we are not aware for which 
constant $c_t>0$ the phase of transition from polynomial to logarithmic behavior occurs. Our 
bounds strongly suggest that the value of this constant seems plausible to be $c_t=1/\log 2$.
\end{note}

\medskip

We close this section with some discussion on the methods used for bounding the variance. If we are interested in 
giving sufficient upper bounds, we may use the Poincar\'{e} inequality \eqref{eq:Poin-ineq} which estimates the variance 
by the $L_2$ average of the Euclidean norm of the gradient of $f$. In principle the latter average is smaller than the Lipschitz 
constant: $\big\| \|\nabla f\|_2 \big\|_{L_2(\gamma_n)} \leq 
\big\| \|\nabla f\|_2 \big\|_{L_\infty (\gamma_n) }=L$. The reader may check that for $2< p<\infty$ and $f=\|\cdot\|_p$ we have:
\begin{align*}
\int_{\mathbb R^n} \big \| \nabla f(x)  \big\|_2^2 \, d\gamma_n(x) = \int_{\mathbb R^n} \frac{\|x\|_{2p-2}^{2p-2}}{\|x\|_p^{2p-2}}\, d\gamma_n(x) \simeq_p \frac{1}{n^{1-\frac{2}{p}}} \ll 1= b^2(B_p^n) \equiv {\rm Lip}(f)^2.
\end{align*} In case $p=\infty$ we have $\big\| \nabla\|x\|_\infty \big\|_2\equiv 1$ a.e., hence:
\begin{align*}
\int_{\mathbb R^n} \big\| \nabla \|x\|_\infty \big\|_2^2 \, d\gamma_n(x) =1 =b(B_\infty^n).
\end{align*} Thus, the Poincar\'{e} inequality also fails to give the sharp upper bound for the variance in this case. 
The recovery of the correct estimate is promised by the different order of magnitude for the $L_1-L_2$ norms 
of the partial derivatives of $x\mapsto \|x\|_\infty$ and
Talagrand's inequality (see \cite{Cha} for the details). The phenomenon that ${\rm Var}\|X\|_{\ell_\infty^n}\simeq 1/\log n$ 
while $\mathbb E \big\| \nabla \|X\|_\infty \big\|_2^2\simeq 1$ is referred to {\it super-concentration} following \cite{Cha}. For recent
results on the related subject see \cite{Tan}.


\section{Gaussian concentration for $\ell_p$ norms}

In this Section we study the Gaussian concentration for the $\ell_p$-norms for $1\leq p\leq \infty$. First we show how 
we may employ the log-Sobolev inequality in order to get concentration results. 

\subsection{An argument via the log-Sobolev inequality}

Note that for the $\ell_p$ norm with $1\leq p\leq 2$ 
the estimate \eqref{eq:2.11} implies:
\begin{align*}
\frac{I_r(\gamma_n,B_p^n)}{I_1(\gamma_n,B_p^n)} \leq \sqrt{1+\frac{C_1r}{k(B_p^n) } } \leq \exp \left( \frac{C_2r}{n} \right),
\end{align*} for all $r\geq 1 $. Therefore, for any $0<\varepsilon <1$ we apply Markov's inequality to get:
\begin{align*}
P(\|X\|_p >(1+\varepsilon) I_1 ) \leq P( \|X\|_p > e^{\varepsilon/2} I_1) \leq e^{-\varepsilon r/2} (I_r/I_1)^r\leq \exp(-\varepsilon r/2 +C_2r^2/n).
\end{align*} Choosing $r= \varepsilon n/(4C_2)$ (as long as $\varepsilon >4C_2/n$) we obtain:
\begin{align*}
P(\|X\|_p >(1+\varepsilon) I_1 ) \leq \exp \left(-\frac{1}{16C_2} \varepsilon^2 n \right).
\end{align*} Taking into account Theorem \ref{thm:stability-moms} and arguing similarly we find:
\begin{align*}
P(\|X\|_p <(1-\varepsilon)I_1 )\leq \exp(-c_2\varepsilon^2n).
\end{align*} Combining those two estimates we arrive at the next concentration result:
\begin{align*}
P\left( \big| \|X\|_p -I_1 \big| >\varepsilon I_1 \right )\leq C_3 \exp(-c_3\varepsilon^2n),
\end{align*} for all $0<\varepsilon <1$. This estimate is sharp, as we will show later, but the same method fails for the $\ell_p$ norm,
when $2<p \leq \infty$, to give the correct concentration estimate. By carefully inspecting the 
proof of the estimates we used before we see that we have bounded the $L_2$ norm of the gradient by the 
$L_\infty$ norm, i.e. the Lipschitz constant. A first attempt to improve the estimates, would be to improve the
bound on that quantity. To this end, we restrict ourselves to the range $2<p<\log n $ and we use the log-Sobolev inequality.
We have the following:

\begin{proposition} \label{prop: stab-ell_p}
Let $2<p< c \log n$. Then, for every $r>0$ we have:
\begin{align*} \displaystyle
\frac{d}{dr}(\log I_r) \leq \frac{C^p}{n} \left( 1+ \frac{r}{k(B_{2p-2}^n)}\right)^{p-1} 
\leq \left\{
\begin{array}{ll}
 C_1^p/n , & 0<r\leq k(B_{2p-2}^n) \\
 \\
\displaystyle  \frac{1}{r} \left( \frac{C_1r}{k(B_p^n)} \right)^p , & k(B_{2p-2}^n) \leq r < k(B_p^n)/C_1
\end{array} \right. ,
\end{align*} while for $0< r < c d(B_p^n)$ we have:
\begin{align*}
- \frac{d}{dr}(\log I_{-r}) \leq \frac{C^p}{n},
\end{align*} where $c,C,C_1>0$ are absolute constants and $I_s\equiv I_s(\gamma_n, B_p^n)$.
\end{proposition}

\noindent {\it Proof.} First we prove the growth condition on the positive moments. Our starting point is the next
estimate:
\begin{align*}
\frac{d}{dr} (\log I_r) =\frac{1}{r^2 I_r^r} {\rm Ent}_{\gamma_n}(\|x\|_p^r) \leq \frac{2}{r^2I_r^r} \mathbb E \left\|  \nabla (\|X\|_p^{r/2})\right\|_2^2 
=\frac{1}{2 I_r^r} \mathbb E\|X\|_{2p-2}^{2p-2} \|X\|_p^{r-2p},
\end{align*} where we have used the log-Sobolev inequality. We distinguish two cases:

\smallskip

\noindent {\it Case 1: $0<r \leq 2p$.} We may write: 
\begin{align*}
\frac{d}{dr} (\log I_r) &\leq  \frac{n}{\mathbb E \|X\|^r_p} 
\mathbb E \frac{|X_1|^{2p-2}}{\|X\|_p^{2p-r}} \leq \frac{n \sigma_{2p-2}^{2p-2}}{ \mathbb E \|X\|^r_{p}} 
\mathbb E \frac{ 1}{\|X\|^{2p-r}_{p} } \leq \frac{n (cp)^{p-1}}{I_r^r(B_p^{n}) I_{-(2p-r)}^{2p-r}(B_p^{n})} 
\leq \frac{n(cp)^p}{I_{-2p}^{2p}(B_p^{n})},
\end{align*} by Proposition \ref{prop:Harris} and H\"{o}lder's inequality.  By Proposition \ref{prop:small-ball} for $0< s < c_1k_{p,n}$ 
we have: $I_{-s}\geq c_2 I_1$. Since, $p < c_1 k_{p,n}$ for $p\lesssim \log n$ we get: $(\log I_r)' \leq C_2^p / n$.

\smallskip

\noindent {\it Case 2: $r> 2p$.} We may write:
\begin{align*}
\frac{d}{dr} (\log I_r) &\leq \frac{1}{2 I_r^r} \mathbb E\|X\|_{2p-2}^{2p-2} \|X\|_p^{r-2p} \leq \frac{I_r^{2p-2}(\gamma_n,B_{2p-2}^n) }{2I_r^{2p}} ,
\end{align*} by H\"{o}lder's inequality. By Lemma \ref{lem:log-sob-moms} we get:
\begin{align*}
\frac{d}{dr} (\log I_r) \leq \frac{I_{2p-2}^{2p-2}(\gamma_n,B_{2p-2}^n) }{2I_p^{2p}} \left(1+\frac{r}{k_{2p-2,n}} \right)^{p-1} 
&=\frac{\sigma_{2p-2}^{2p-2}/\sigma_p^{2p}}{2n} \left(1+\frac{r}{k_{2p-2,n} }\right)^{p-1} \\
&\leq \frac{C_3^p}{n} \left(1+\frac{r}{k_{2p-2,n} }\right)^{p-1},
\end{align*} for some absolute constant $C_3>0$. 

\smallskip

Now we turn to providing bounds for the negative moments. Here the argument is simpler. Using the log-Sobolev inequality
again and Proposition \ref{prop:Harris} we have:
\begin{align*}
\frac{d}{dr} (\log I_{-r}) & \geq - \frac{1}{2I_{-r}^{-r}} \mathbb E \|X\|_{2p-2}^{2p-2} \|X\|_p^{-r-2p} 
\geq - \frac{1}{2I_{-r}^{-r}} \mathbb E \|X\|_{2p-2}^{2p-2} \mathbb E\|X\|_p^{-r-2p} \\
&= - \frac{1}{2I_{-r}^{-r}} \mathbb E \|X\|_{2p-2}^{2p-2} I_{-r-2p}^{-r-2p} \geq  -C_2^p \frac{ \sigma_{2p-2}^{2p-2} n}{I_1^{2p}} \geq -C_3^p/n ,
\end{align*} for $r\leq c_4 d(B_p^n)$, where in the last step we have used Lemma \ref{lem:reduc-neg-moms}. The result easily 
follows.  \prend

\medskip

We are ready to prove the next concentration inequality. Note that the dependence we get on $\varepsilon$ is
better than the one we get if we employ \eqref{eq:2.14}.

\begin{proposition} \label{prop: weak-conc}
Let $4\leq p < c_0 \log n$. Then, one has:
\begin{align*}
P \left( \big | \|X\|_p- \mathbb E\|X\|_p \big| >\varepsilon \mathbb E\|X\|_p \right) \leq C_1\exp \left(-c_1 \varepsilon^{1+\frac{1}{p}} k(B_p^n) \right),
\end{align*} for all $0<\varepsilon<1$. Moreover, we have:
\begin{align*}
P \left( \|X\|_p \leq (1-\varepsilon) \mathbb E \|X\|_p \right) \leq C_2 \exp \left( -c_2\varepsilon k(B_p^n) \right),
\end{align*} for $0<\varepsilon <1$.
\end{proposition}

\noindent {\it Proof.} Let $4\leq p\leq c\log n$, where $c>0$ is the constant from Proposition \ref{prop: stab-ell_p}. 
Then, for each $0<\varepsilon<1$ using Markov's inequality we may write:
\begin{align} \label{eq:4.8}
P(\|X\|_p > (1+\varepsilon) I_0 )\leq e^{-\varepsilon r/2} \exp( r \log (I_r/I_0)) =\exp\left[ -r\left(\frac{\varepsilon}{2}- \log(I_r/I_0) \right) \right],
\end{align} for all $r>0$. Using Proposition \ref{prop: stab-ell_p} we obtain:
\begin{align*}
\log (I_r/I_0) \leq \frac{C^p}{n}\int_0^r \left( 1+\frac{s}{k_{2p-2,n}} \right)^{p-1} \, ds < \frac{C^p k_{2p-2,n} }{pn } \left( 1+\frac{r}{k_{2p-2,n}} \right)^p 
<  \frac{(2C)^p k_{2p-2,n} }{pn } \left(\frac{r}{k_{2p-2,n}} \right)^p ,
\end{align*} for $r> k_{2p-2,n}$. Therefore, \eqref{eq:4.8} becomes:
\begin{align*}
P(\|X\|_p > (1+\varepsilon) I_0 )\leq \exp \left(-\frac{\varepsilon r}{2} +\frac{(2C)^p }{pn k_{2p-2,n}^{p-1}} r^{p+1} \right),
\end{align*} for $r>k_{2p-2,n}$.
 Minimizing the right-hand side with respect to $r$, we find that $r_{\rm min} =r_0$ satisfies:
 \begin{align} \label{eq:4.10}
 \frac{(2C)^p }{pn k_{2p-2,n}^{p-1}} (p+1) r_0^p -\frac{\varepsilon}{2} =0 \Longrightarrow r_0\simeq \varepsilon^{1/p} k_{p,n},
\end{align} and in order for this value to be admissible we ought to have $r_0>k_{2p-2,n}$. Hence, the value $r_0$ is 
admissible if $\varepsilon$ satisfies:
\begin{align*}
r_0>k_{2p-2,n} \Longleftrightarrow (2C)^{-p} \frac{\varepsilon n}{2} \frac{p}{p+1} > k_{2p-2,n} \Longleftrightarrow \varepsilon >(2C)^p \frac{2(p+1)}{pn} k_{2p-2,n}.
\end{align*} Note that Proposition \ref{prop:mean-ell-p} implies that:
\begin{align} \label{eq: 4.11}
k_{q,n} \leq c_2qn^{2/q}, \quad \forall \, 2\leq q\leq \log n.
\end{align}
Since $p\geq 4$ it suffices to have $\varepsilon > (2C)^p 8c_2 p n^{-\frac{p-2}{p-1}}$ or equivalently
to have $\varepsilon > (16ec_2C)^p n^{-\frac{p-2}{p-1}}$. 

First consider the case $k_{p,n}^{-\frac{p}{p+1}} <\varepsilon <1$. In this case the above restriction is satisfied as long as
$p\leq c_3\log n$ for some sufficiently small absolute constant $c_3>0$. Indeed one needs to check that:
$k_{p,n}^{-\frac{p}{p+1}} > (16ec_2C)^p n^{-\frac{p-2}{p-1}}$ and by taking into account \eqref{eq: 4.11} again it suffices to have 
$\frac{ n^{\frac{p-2}{p-1}} }{(c_2p n^{2/p})^{\frac{p}{p+1}}} > (16ec_2C)^p$ or it is enough to have:
\begin{align*}
n^{\frac{p^2-3p}{p^2-1}} > (16e^2c_2^2 C)^p=e^{p/c_4}.
\end{align*}
Thus, if $c_0:= \min \{c_3,c_4/4\}>0$ we have all the requirements so that we may conclude:
\begin{align*}
P(\|X\|_p > (1+\varepsilon) I_0 )\leq \exp \left(-\frac{\varepsilon r_0}{2} +\frac{(2C)^p }{pn k_{2p-2,n}^{p-1}} r_0^{p+1} \right) 
& \stackrel{ \eqref{eq:4.10}}=\exp \left( -\frac{\varepsilon}{2}r_0 +\frac{\varepsilon r_0}{2(p+1)} \right) \\
&=\exp \left( -\frac{p}{2(p+1)} \varepsilon r_0 \right) \\
& \leq \exp \left(-c \varepsilon^{1+\frac{1}{p}} k_{p,n} \right),
\end{align*} for all $4\leq p\leq c_0 \log n$ and for all $k_{p,n}^{-\frac{p}{p+1}} <\varepsilon<1$. By adjusting the constants we get:
\begin{align*}
P(\|X\|_p > (1+\varepsilon) I_0 )\leq C' \exp \left(-c \varepsilon^{1+\frac{1}{p}} k_{p,n} \right),
\end{align*} for the whole range $0<\varepsilon <1$ and for $4\leq p \leq c_0 \log n$. 

\medskip
 
\noindent Now we turn to bounding the probability $P(\|X\| \leq (1-\varepsilon)I_0)$. 
Proposition \ref{prop: stab-ell_p} shows that $(\log I_{-r})' \geq -C^p/n$ for $0<r\leq c_1k_p$. Hence, we get:
\begin{align*}
P(\|X\|_p \leq (1-\varepsilon)I_0) \leq P( \|X\|_ p\leq e^{-\varepsilon} I_0) \leq e^{-\varepsilon r} \left( \frac{I_0}{I_{-r}} \right)^r \leq \exp(-\varepsilon r+ r^2C^p/n),
\end{align*} for all $0<r<c_1k_{p,n}$, where we have used the bound:
\begin{align*}
\log(I_0/I_{-r}) = -\int_0^r (\log I_{-s})' \, ds \leq \frac{C^p}{n}r,
\end{align*} for $0<r<c_1k_{p,n}$. Finally, choosing $r\simeq k_{p,n}$ we see 
that $C^pk_{p,n}^2/n< (2eC)^p n^{4/p-1} \leq C'$ as long as $4\leq p\leq c_1' \log n$, hence we conclude:
\begin{align*}
P(\|X\|_p \leq (1-\varepsilon)I_0) \leq C' \exp(-c' \varepsilon k_{p,n}),
\end{align*} for $0<\varepsilon <1$. \prend

\medskip

Although this concentration result improves upon the one we get by just using \eqref{eq:2.14}, it is still suboptimal. It turns out that
although the $L_2$ average of the Euclidean norm of the gradient is the proper quantity to be estimated for the concentration result,
yet it should not be used in order to bound the growth of the high moments of the norm, in this range of $p$.

\subsection{Estimating centered moments}

In this paragraph we study centered moments of the $\ell_p$ norm, i.e. $(\mathbb E \left| \|X\|_p- \|Y\|_p \right|^r)^{1/r}$. 
For this end we distinguish three cases: (a) $1\leq p\leq 2$, (b) $2<p < c_0\log n$ and (c) $c_0\log n\leq p\leq \infty$, where $c_0>0$ is sufficiently small absolute constant. 
While in the first two cases we estimate directly the centered moments in terms of $n,p,r$, in the last we have to argue
differently and study the almost constant behavior of the noncentered moments. This is because when $p$ grows along with $n$ 
the estimates collapse. To overcome this obstacle we use Talagrand's $L_1-L_2$ bound.

\subsubsection{The case $1\leq p\leq 2$}

In this subsection we sketch the proof of the next theorem:

\begin{theorem} \label{thm:conc-1<p<2}
Let $1\leq p\leq 2$. Then, one has:
\begin{align}\label{eq:4.1}
c_1\exp(-C_1\varepsilon^2n)\leq P\left( \left| \|X\|_p -\mathbb E\|X\|_p\right| > \varepsilon \mathbb E\|X\|_p \right)
\leq C_2\exp(-c_2\varepsilon^2 n),
\end{align} for $0 < \varepsilon <1$, where $C_1,c_1,C_2,c_2>0$ are absolute constants.
\end{theorem}

\noindent {\it Proof (Sketch).}  The rightmost inequality follows by the Gaussian concentration 
inequality \eqref{eq:2.14}, Proposition \ref{prop:mean-ell-p}
and the fact that ${\rm Lip}(\|\cdot\|_p)=b(B_p^n)=n^{1/p-1/2}$ for $1\leq p\leq 2$. Now we focus on 
the left-hand side inequality. We have the next:

\begin{proposition} Let $1\leq p\leq 2$. Then, we have:
\begin{align*}
\left( \mathbb E \left| \|X\|_p-\mathbb E\|X\|_p \right|^r \right)^{1/r} \simeq \sqrt{\frac{r}{n}} \mathbb E\|X\|_p,
\end{align*} for all $r\geq 1$.
\end{proposition}

\noindent {\it Proof.} Indeed; the estimate 
\begin{align} \label{eq:4.3} \left( \mathbb E \left| \|X\|_p-\mathbb E\|X\|_p \right|^r \right)^{1/r} \leq C_3 \sqrt{\frac{r}{n}} \mathbb E\|X\|_p , \quad r\geq 1
\end{align} is well known and follows by integration by parts combined with the right-hand side estimate in \eqref{eq:4.1}.
For the estimate 
\begin{align} \label{eq:4.4}
\left( \mathbb E \left| \|X\|_p-\mathbb E\|X\|_p \right|^r \right)^{1/r} \geq c_3 \sqrt{\frac{r}{n}} \mathbb E\|X\|_p
\end{align}
we may apply the triangle inequality, Lemma \ref{lem:2-sided-ineq} and finally the Cauchy-Schwarz inequality to write:
\begin{align*}
2\left( \mathbb E \left| \|X\|_p-\mathbb E\|X\|_p \right|^r \right)^{1/r} \geq \left( \mathbb E \left| \|X\|_p-\|Y\|_p \right|^r \right)^{1/r}
 \geq \frac{1}{2p} \frac{\left(\mathbb E\big|\|X\|_p^p -\|Y\|_p^p \big|^{r/2}\right)^{2/r}}{\left( \mathbb E \|X\|_p^{r(p-1)}\right)^{1/r} }.
\end{align*} Note that \eqref{eq:4.3} already implies
$\left(\mathbb E\|X\|_p^s\right)^{1/s} \leq 2C_3 \mathbb E\|X\|_p\simeq n^{1/p}$ for all $1\leq s\leq n$. Moreover, we have: 
\begin{align*}
\left(\mathbb E \left| \|X\|_p^p -\|Y\|_p^p \right|^s\right)^{1/s} 
\geq \mathbb E\big| |X_1|^p-|Y_1|^p \big| \cdot \left( \mathbb E_\varepsilon \left| \sum_{i=1}^n \varepsilon_i \right|^s \right)^{1/s}\simeq \sqrt{sn},
\end{align*} where we have used the facts that the joint distribution of $(\varepsilon_i | |X_i|^p-|Y_i|^p |)_i$ is the same as 
$(|X_i|^p-|Y_i|^p)_i$, Jensen's inequality and at the last step, that 
$\left( \mathbb E_\varepsilon \left| \sum_{i=1}^n \varepsilon_i \right|^s \right)^{1/s} \simeq \sqrt{sn}$ for $1\leq s\leq n$ 
(see e.g. \cite{M}). 
Putting them all together we see:
\begin{align*}
\left( \mathbb E \left| \|X\|_p-\mathbb E\|X\|_p \right|^r \right)^{1/r} \geq c_4 \frac{\sqrt{rn}}{n^{1-1/p}} 
\simeq \sqrt{ \frac{r}{n}} \mathbb E\|X\|_p,
\end{align*} which completes the proof. \prend 

\smallskip

Now we turn to proving the lower bound in the probabilistic estimate \eqref{eq:4.1}:
For every $n^{-1/2} <\varepsilon < 2c_3$ consider $r\in [1,n]$ so that $\varepsilon= 2c_3\sqrt{r/n}$ to write:
\begin{align*}
P \left( \big| \|X\|_p -\mathbb E \|X\|_p \big| >\varepsilon \mathbb E\|X\|_p \right) & \geq 
P \left( \big| \|X\|_p-\mathbb E\|X\|_p \big| \geq \frac{1}{2} \left( \mathbb E\big| \|X\|_p-\mathbb E\|X\|_p \big|^r \right)^{1/r} \right) \\
&= P(\zeta \geq 2^{-r} \mathbb E\zeta) \geq (1-2^{-r})^2 \frac{(\mathbb E\zeta)^2}{\mathbb E\zeta^2},
\end{align*} by Lemma \ref{lem: PZ-ineq}, where $\zeta := \big| \|X\|_p-\mathbb E\|X\|_p \big|^r$. Employing the estimates \eqref{eq:4.3} and 
\eqref{eq:4.4} we conclude:
\begin{align*}
P \left( \big| \|X\|_p -\mathbb E \|X\|_p \big| >\varepsilon \mathbb E\|X\|_p \right)  \geq c_5 e^{-C_5r},
\end{align*} as required. \prend

\subsubsection{The case $2<p \leq c_0 \log n$}
It is clear from the argument of the previous paragraph that in order
to obtain sharp concentration inequalities it is enough to get sharp estimates for the centered moments:
$\left( \mathbb E \left| \|X\|_p-\|Y\|_p \right|^r \right)^{1/r}.$ 
In view of Lemma \ref{lem:2-sided-ineq} it is also obvious that estimates for the centered moments 
$\left( \mathbb E \big| \|X\|_p^p- \|Y\|_p^p \big|^r \right)^{1/r}$ will provide estimates for the moments 
$\left( \mathbb E \left| \|X\|_p-\|Y\|_p \right|^r \right)^{1/r}$. 
Note that in order to estimate the centered moments from above one may also employ Theorem \ref{thm:Pis-ineq} in the form of an 
$(r,r)$-Poincar\'{e} inequality \eqref{eq: 2.9}. We use this method in the next Section in order to derive the optimal dependence on 
$\varepsilon$ in the critical dimension of randomized Dvoretzky. Here we shall prove the next result (see \cite{Naor} for a different approach):

\begin{proposition} \label{prop:ctrd-p-pwr}
Let $2< p <\infty$. Then, we have:
\begin{align}
\left( \mathbb E \big| \|X\|_p^p- \|Y\|_p^p \big|^r \right)^{1/r} \simeq 
\sigma_p^p \max \left\{ 2^{p/2} ( r n )^{1/2} , \, r^{p/2} n^{1/r} \right\} ,
\end{align} for all $r\geq 2$.
\end{proposition} 

\noindent {\it Proof.} Note that if $X=(X_1,\ldots,X_n)$
is a Gaussian random vector and $Y$ an independent copy of it, the variables $\xi_i:= |X_i|^p-|Y_i|^p$ are i.i.d. and the 
functions $t\mapsto P(|\xi_i| >t)$ are log-convex on $(0,\infty)$ by Lemma \ref{lem:gauss^p}. Then we may apply the main 
result from \cite{HMSO} to get:
\begin{align*} 
\left( \mathbb E \big| \|X\|_p^p -\|Y\|_p^p \big|^r\right)^{1/r} \equiv \left\| \sum_{i=1}^n \xi_i \right\|_r &\simeq 
\left(\sum_{i=1}^n\|\xi_i\|_r^r \right)^{1/r} + \sqrt{r} \left( \sum_{i=1}^n \| \xi_i\|_2^2 \right)^{1/2} \\
&\simeq n^{1/r} \|\xi_1\|_r +\sqrt{rn} \|\xi_1\|_2 \\
&\simeq n^{1/r} r^{p/2} \sigma_p^p + \sqrt{rn} 2^{p/2} \sigma_p^p,
\end{align*} where we have used Lemma \ref{lem:gauss^p} again. The proof is complete. \prend

\medskip

Now we are ready to prove the following:

\begin{theorem} \label{thm:conc-p>2-I-p}
Let $n>C$ and let $2<p \leq c_0 \log n$. Then, we have:
\begin{align*}
P\left( \left| \|X\|_p - \mathbb E \|X\|_p \right| > \varepsilon \mathbb E\|X\|_p \right)\leq C \exp \left(-c \min \left\{ \frac{\varepsilon^2 p^2 n}{2^p}, (\varepsilon n)^{2/p} \right\} \right),
\end{align*} for all $0<\varepsilon<1/p$, where $C,c,c_0>0$ are absolute constants.
\end{theorem}

\noindent {\it Proof.} Define $\alpha(n,p,r):= \max\{ 2^{p/2} (rn)^{1/2}, r^{p/2}n^{1/r} \}, \; r\geq 2$. Note that for fixed
$n,p$ the map $r \mapsto \alpha(n,p,r)$ is strictly increasing with inverse $A(n,p,s) \simeq \min\{ \frac{s^2}{2^p n} , s^{2/p}\}$. 
Then, Proposition \ref{prop:ctrd-p-pwr} shows that:
\begin{align} \label{eq:equiv-p-moms}
c_1 \sigma_p^p \alpha(n,p,r) \leq \left( \mathbb E \big| \|X\|_p^p- \mathbb E\|X\|_p^p \big|^r \right)^{1/r} \leq C_1
\sigma_p^p \alpha(n,p,r),
\end{align} for all $r\geq 2$. Applying Markov's inequality we get:
\begin{align*}
P\left( \left| \|X\|_p^p -\mathbb E \|X\|_p^p \right| > t \mathbb E\|X\|_p^p \right) \leq \left( \frac{C_1 \alpha(n,p,r) }{t n} \right)^r 
&= \exp \left( - A(n,p, etn/C_1)\right) \\
&\leq \exp \left(-c_2 \min\left\{ \frac{t^2 n}{2^p}, (tn)^{2/p} \right\} \right) ,
\end{align*} provided that $etn/C_1 > \alpha(n,p,2)=2^{p/2}n^{1/2}$. It follows that
\begin{align*}
P\left( \left| \|X\|_p^p -\mathbb E \|X\|_p^p \right| > t \mathbb E\|X\|_p^p \right) \leq e^2 \exp \left(-c_2 \min\left\{ \frac{t^2 n}{2^p}, (tn)^{2/p} \right\} \right),\end{align*} for all $t>0$. Now fix $0< \varepsilon<1/p$. Then, we may write:
\begin{align*}
P \left( \|X\|_p < (1-\varepsilon) ( \mathbb E\|X\|_p^p)^{1/p} \right) \leq P\left( \|X\|_p^p < \left(1- \frac{\varepsilon p}{2} \right) \mathbb E\|X\|_p^p \right) 
\leq e^2 \exp \left(-c_3 \min \left\{ \frac{\varepsilon^2 p^2 n}{2^p}, (\varepsilon n)^{2/p} \right\}\right),
\end{align*} by the previous estimate. Arguing similarly we show the upper estimate. Thus,
\begin{align*}
P\left( \left| \|X\|_p - (\mathbb E \|X\|_p^p)^{1/p} \right| > \varepsilon (\mathbb E\|X\|_p^p)^{1/p} \right)\leq e^2 \exp \left(-c_4 \min \left\{ \frac{\varepsilon^2 p^2 n}{2^p}, (\varepsilon n)^{2/p} \right\}\right),
\end{align*} for all $0<\varepsilon<1/p$. The result follows.  \prend

\medskip

\begin{remark}\rm  For fixed $2<p< \infty$ the estimate can be reversed. The argument is similar to that in Theorem \ref{thm:conc-1<p<2}.
Let $\frac{ 2^{p/2} }{2c_1 n^{1/2} } <t <1$ and choose $r\geq 2$ with $\alpha(n,p,r) = 2c_1 n t$. Then, in view of the lower estimate 
in \eqref{eq:equiv-p-moms} we get:
\begin{align*}
P \left(\big| \|X\|_p^p-\mathbb E\|X\|_p^p \big| > t \mathbb E\|X\|_p^p \right) &\geq 
P \left(\big| \|X\|_p^p-\mathbb E\|X\|_p^p \big| > \frac{1}{2}  \left( \mathbb E \left| \|X\|_p^p -\mathbb E \|X\|_p^p \right|^r \right)^{1/r} \right) \\
&\geq (1-2^{-r})^2 \frac{\left( \mathbb E \big| \|X\|_p^p-\mathbb E\|X\|_p^p \big|^r \right)^2}{\mathbb E \big| \|X\|_p^p-\mathbb E\|X\|_p^p \big|^{2r} } \\
& \geq \frac{1}{4} e^{-c_5 p r} \geq \frac{1}{4} \exp\left(-c_6 p A(n, p, 2c_1nt) \right).
\end{align*}  It follows, as before, that:
\begin{align*}
P \left(\big| \|X\|_p-(\mathbb E\|X\|_p^p)^{1/p} \big| > \varepsilon (\mathbb E\|X\|_p^p)^{1/p} \right) \geq c_7\exp\left(-c_8 p A(n,p, 2c_1np\varepsilon) \right),
\end{align*} for $ \frac{2^{p/2}}{2c_1p n^{1/2}} < \varepsilon <1/p$. Next recall that Proposition \ref{prop: stab-ell_p} implies $I_p/I_1 \leq 1+ \frac{2^{p/2}}{4c_1p n^{1/2}}$ for $p\leq c_0\log n$, where $I_r\equiv I_r(\gamma_n,B_p^n)$. Thus, we may replace $I_p$ by $I_1$ in the above concentration estimate. This yields the following double estimate:

\begin{proposition} \label{prop:conc-p>2-u-l-b}For all sufficiently large $n$ and for $2<p<c_0\log n$ one has:
\begin{align}
c\exp\left( -C p \theta(n,p,\varepsilon) \right) \leq P\left( \left| \|X\|_p - \mathbb E \|X\|_p \right| > \varepsilon \mathbb E\|X\|_p \right) 
\leq C \exp \left(-c \theta(n,p,\varepsilon) \right),
\end{align} for all $0<\varepsilon<1/p$, where $\theta(n,p,\varepsilon) := \min \left\{ \frac{p^2 \varepsilon^2 n}{2^p}, (\varepsilon n)^{2/p}\right\}$ and
$C, c, c_0>0$ are absolute constants.
\end{proposition}

\end{remark}

\noindent {\it Note.} Let us mention that the extra $p$ on the exponent in the lower estimate can be removed if we restrict the range to
$p^{-1}2^{p/2}n^{-1/2} \lesssim \varepsilon \lesssim p^{-1}2^{p/2} n^{-\frac{p-2}{2(p-1)} }$.

\subsubsection{The case $c_0\log n< p\leq \infty$}

In this Subsection we deal with the large values of $p$ in terms of the dimension, namely for $p\gtrsim \log n$. We have the following:

\begin{theorem}\label{thm:stability-r-means}
Let $4< p\leq \infty$. Then, for any $0< r <s \leq c_1\sqrt{k_{p,n} \log n}$ we have:
\begin{align*}
\frac{I_s(\gamma_n,B_p^n)}{I_r(\gamma_n,B_p^n)} \leq \exp \left( \frac{c_2(2s-r)}{k_{p,n} \log n} \right), \quad
\frac{I_{-s} (\gamma_n,B_p^n)}{I_{-r}(\gamma_n,B_p^n)} \geq \exp \left( -\frac{c_2(2s-r)}{k_{p,n} \log n} \right),
\end{align*} where $c_1,c_2>0$ are absolute constants.
\end{theorem}

\noindent {\it Proof.} Set $I_s\equiv I_s(\gamma_n, B_p^n)$. If $a= a_i := \|\partial_i f\|_{L_1(\gamma_n)} $ we get:
\begin{align} \label{eq:L-1-conc}
a_i = \frac{|r|}{n} \int_{\mathbb R^n} \|x\|_{p-1}^{p-1} \|x\|_p^{r-p}\, d\gamma_n(x) \leq  \frac{|r|}{n^{1/q}}I_{r-1}^{r-1},
\end{align} where we have used \eqref{eq: Holder - p-norms}. Similarly, for $A=A_i:= \|\partial_i f\|_{L_2(\gamma_n)}$ we have that:
\begin{align} \label{eq:L-2-conc}
\frac{|r|}{n^{1/q}} I_{2r-2}^{r-1} \leq A= \frac{|r|}{n^{1/2}} \left( \int_{\mathbb R^n} \|x\|_p^{2r-2p} \|x\|_{2p-2}^{2p-2} \, d\gamma_n(x) \right)^{1/2} \leq 
\frac{|r|}{n^{1/2}} I_{2r-2}^{r-1}.
\end{align} 
We apply Theorem \ref{thm:Talagrand bd} for $f(x):=\|x\|_p^r, \, r\neq 0$ to obtain:
\begin{align*}
{\rm Var}_{\gamma_n}(f) \leq C_1 n \frac{A^2}{ 1+\log(A/a)} .
\end{align*}  The function $t \mapsto \frac{t^2}{1+\log (t/a)}, \; t>a$ is increasing, thus  \eqref{eq:L-1-conc} and \eqref{eq:L-2-conc} imply that: 
\begin{align} \label{eq: 4.35}
I_{2r}^{2r}-I_r^{2r} = {\rm Var}_{\gamma_n} (f) \leq C_1 r^2 \frac{I_{2r-2}^{2r-2}}{1+ \log \left(n^{1/q-1/2} \left(\frac{I_{2r-2}}{I_{r-1}} \right)^{r-1} \right) }
\leq C_2 r^2 \frac{I_{2r-2}^{2r-2}}{\log n},
\end{align} for all $r\neq 0$, since $1\leq q<4/3$ and $\log (I_{2r-2}/I_{r-1}) \geq 0$.

\smallskip

\noindent {\it Claim.} For $r > - k_{p,n}, \; r\neq 0$ we have: 
\begin{align*}I_{2r-2}^{2r-2} \leq C_3 I_{2r}^{2r} / k_{p,n}.
\end{align*}
We distinguish three cases:
\begin{itemize}
\item For $0<r<1$ we have: $I_{2r-2}^{2r-2} =\frac{I_{2r-2}^{2r}}{I_{2r-2}^2}\leq  c_2' I_1^{-2} I_{2r-2}^{2r} \leq \frac{c_2'}{k_{p,n} } I_{2r}^{2r}$.
\item For $r\geq 1$ we may write: $I_{2r-2}^{2r-2} = \frac{I_{2r-2}^{2r} }{ I_{2r-2}^2 } \leq  c_3 \frac{ I_{2r}^{2r} }{I_1^2} = \frac{c_3}{k_{p,n}} I_{2r}^{2r}$, 
since $I_1\simeq I_0$. 
\item Finally, for $- k_{p,n} < r<0$ we have: $I_{2r-2}^{2r-2} \leq \frac{c_4}{I_1^2} I_{2r}^{2r} =\frac{c_4}{k_{p,n} } I_{2r}^{2r}$, by 
Lemma \ref{lem:reduc-neg-moms}.
\end{itemize} Thus, \eqref{eq: 4.35} yields:
\begin{align} \label{eq:recursive-ineq}
I_{2r}^{2r} -I_r^{2r} \leq Cr^2 \frac{I_{2r}^{2r}}{k_{p,n} \log n}.
\end{align} for $r> - k_p, \; r\neq 0$. We only prove the stability for the positive moments (the negative moments are treated similarly):
As long as $0 < r < \sqrt{k_{p,n} \log n /C}$ we may write
\begin{align*}
I_{2r}^{2r} \leq \left( 1+\frac{Cr^2}{k_{p,n} \log n} \right) I_r^{2r}.
\end{align*} Iterating the last one we find:
\begin{align*}
\frac{I_{2^mr}}{I_r} \leq \exp \left( C\sum_{j=0}^{m-1} \frac{2^j r}{k_{p,n} \log n}\right) \leq \exp\left( \frac{C(2^mr-r)}{k_{p,n} \log n}\right),
\end{align*}  for $m=1,2,\ldots$ as long as $2^mr \leq \sqrt{k_{p,n} \log n/C}$. 
The result follows. \prend

\medskip

The next corollary is immediate:

\begin{corollary} \label{cor:conc-large-p}
Let $ c_0\log n<p\leq \infty$. Then, one has:
\begin{align*}
P \left( \big| \|X\|_p -\mathbb E\|X\|_p \big| >\varepsilon \mathbb E\|X\|_p \right) \leq C \exp \left(-c \varepsilon \log n \right),
\end{align*} for all $\varepsilon \in (0,1)$, where $C,c,c_0>0$ are absolute constants.
\end{corollary}

\noindent {\it Proof.} Let $K:=k_{p,n} \log n$. Using Markov's inequality and Theorem \ref{thm:stability-r-means} we may write:
\begin{align*}
P \left( \|X\|_p\geq (1+\varepsilon) I_1 \right)\leq P \left( \|X\|_p \geq e^{\varepsilon/2} I_1 \right)\leq e^{-\varepsilon r/2} \left(\frac{I_r}{I_1}\right)^r 
\leq \exp(-\varepsilon r/2+c_2r^2/K),
\end{align*} for all $0<r<c_1\sqrt{K}$. The choice $r\simeq \sqrt{K}$ yields the one-sided estimate:
\begin{align*}
P \left( \|X\|_p >(1+\varepsilon) I_1 \right)\leq C_1\exp \left( -c_1' \varepsilon \sqrt{K} \right).
\end{align*} Working similarly with the probability $P( \|X\|_p < (1-\varepsilon) I_1)$ and taking into account the fact that 
$k(B_p^n) \simeq  \log n$ for $p \gtrsim \log n$, we conclude the asserted estimate. \prend

\medskip

Summarizing the results of this paragraph (by taking into account Theorem \ref{thm:conc-1<p<2}, Theorem \ref{thm:conc-p>2-I-p} and Proposition \ref{prop: weak-conc} and the variance estimates from Section 3) we may have a concentration inequality which interpolates between the concentration estimates for fixed $p\geq 1$ and $p=\infty$:

\begin{theorem} \label{thm:conc-full} For all large enough $n$ and for any $1\leq p\leq \infty$ one has:
\begin{align*}
P\left( \big | \|X\|_p-\mathbb E\|X\|_p \big| > \varepsilon \mathbb E\|X\|_p \right) \leq C_1\exp ( -c_1\beta(n,p,\varepsilon) ),
\end{align*} for every $0<\varepsilon <1$, where $\beta(n,p,\varepsilon)$ is defined as follows:
\begin{align*}
\beta(n,p,\varepsilon) = \left\{ \begin{array}{lll}
 \varepsilon^2 n, & 1\leq p\leq 2 \\
 \max \left\{  \min \left\{ p^2 2^{-p} \varepsilon^2n , (\varepsilon n)^{2/p} \right \} , \varepsilon pn^{2/p} \right\}, & 2< p\leq c_0 \log n \\
 \varepsilon \log n, & p> c_0 \log n
 \end{array}\right. ,
\end{align*}
where $c_0\in (0,1)$ and $C_1,c_1>0$ are suitable absolute constants. Furthermore, for $p\leq c_0 \log n$ we have the estimate:
\begin{align*}
P\left( \big | \|X\|_p-\mathbb E\|X\|_p \big| > \varepsilon \mathbb E\|X\|_p \right) \leq 
\exp \left( -\log\left(1 +c_1\frac{p^2}{2^p} \varepsilon^2n \right) \right),
\end{align*} for every $\varepsilon\in (0,1)$.
\end{theorem}

\section{The critical dimension in random Dvoretzky for $\ell_p^n$}

In this paragraph we study the critical dimension $k(n,p,\varepsilon)$ (and in particular the dependence on $\varepsilon$) in the random version of Dvoretzky's theorem for $\ell_p^n$ spaces.
Our method is inspired by Schechtman's approach in \cite{Sch1}. The key point is a distributional inequality for rectangular 
matrices with independent standard Gaussian entries. In \cite{Sch1} it is proved that, if $G=(g_{ij})_{i,j=1}^{n,k}$ is a Gaussian 
matrix then the process $(\|Gx\|)_{x\in S^{k-1}}$ is sub-Gaussian with constant $b=\max_{\theta\in S^{n-1}}\|\theta\|$.
The proof of \cite[Lemma]{Sch1} is based on an orthogonal splitting, combined with a conditioning argument and 
inequality \eqref{eq:2.14}. 

Here we use similar ideas to prove a functional inequality which generalizes \cite[Lemma]{Sch1}. Once again, the 
advantage of this new inequality is that it involves $\| \nabla f\|_2$ instead of the Lipschitz constant of $f$.
Our result reads as follows:

\begin{theorem}\label{thm:Sch-funct-ineq}
Let $a,b\in S^{k-1}$ and $G=(g_{ij})_{i,j=1}^{n,k}$ be random matrix with standard i.i.d. Gaussian entries. 
If $f:\mathbb R^n\to \mathbb R$ is $C^1$-smooth, then we have:
\begin{align*}
\left( \mathbb E \big| f(Ga)-f(Gb) \big|^r \right)^{1/r} \leq \pi \sigma_r \|a-b\|_2 \left(\mathbb E \|\nabla f(W)\|_2^r \right)^{1/r},
\end{align*} for all $r\geq 1$, where $W\sim N(\mathbf 0, I_n)$. 
 \end{theorem} 

\noindent {\it Proof.}  Fix $a,b\in S^{k-1}$ and assume without loss of generality that $a\neq \pm b$. Define $p:=\frac{a+b}{2}$
and note that since $\|a\|_2=\|b\|_2$ the vector $u:=a-p$ is perpendicular to $p$. 
If we set $X:=G(u)$ and $Z:=G(p)$ then $X, Z$ are independent Gaussian random vectors in $\mathbb R^n$ with
$X \sim N({\bf 0}, \|u\|_2^2I_n)$, $Z \sim N({\bf 0}, \|p\|_2^2I_n)$ and $G(a)=Z+X$ while $G(b)=Z-X$. Thus, we may write:
\begin{align*}
\mathbb E \left| f(Ga) -f(Gb) \right|^r =\mathbb E_Z\mathbb E_X \left| f(Z+X) -f(Z-X) \right|^r.
\end{align*}  For $x, z\in \mathbb R^n$ we define $F(x,z):=f(z+x)-f(z-x)$. Note that for fixed $z$ we have $\mathbb E_X F(X,z)=0$ 
since, $X$ is symmetric random vector. Applying Theorem \ref{thm:Pis-ineq} for $\phi(t)=|t|^r,\; r\geq 1$ and $x\mapsto F(x,z)$ instead of $f$ we derive:
\begin{align*}
\mathbb E|F(X,z)| =
\mathbb E_X \left| f(z+X) -f(z-X) \right|^r  &\leq \left(\frac{\pi}{2}\right)^r \mathbb E_{X,Y} \left| \langle \nabla f(z+X) ,Y\rangle + \langle \nabla f(z-X),Y\rangle \right|^r \\
&\leq \pi^r \mathbb E_{X,Y} \left| \langle \nabla f(z+X), Y\rangle \right|^r\\
&= \pi^r \|a-b\|_2^r \sigma_r^r \, \mathbb E_X \left \|\nabla f(z+X) \right\|_2^r.
\end{align*} Moreover, note that $W:=X+Z\sim N({\bf 0},I_n)$, thus we get:
\begin{align*}
\mathbb E \left| f(Ga) -f(Gb) \right|^r = \mathbb E|F(X,Z)| &\leq  \pi^r \|a-b\|_2^r \sigma_r^r \, \mathbb E_{X,Z} \left \|\nabla f(Z+X) \right\|_2^r,
\end{align*} as required. \prend

\smallskip

\begin{remarks}\rm 1. If we assume that $f$ is $L$-Lipschitz and applying Markov's inequality we may
conclude the more general form of \cite[Lemma]{Sch1}:
\begin{align*} 
{\rm Prob} \left( \big| f(G(a))-f(G(b)) \big| > t \right) \leq 2 \exp\left(-\frac{2}{\pi^2} \frac{t^2}{L^2 \|a-b\|_2^2}\right), \quad t>0.
\end{align*}

\smallskip

\noindent 2. The same proof provides the following variant of Theorem \ref{thm:Pis-ineq} which we state here for future reference:
\begin{theorem}
Let $\phi:\mathbb R\to \mathbb R$ be convex function and let $f:\mathbb R^n\to \mathbb R$ be $C^1$-smooth. If 
$G=(g_{ij})_{i,j=1}^{n,k}$ is Gaussian matrix and $a,b\in S^{k-1}$, then we have:
\begin{align*} 
\mathbb E \phi\left( f(Ga)-f(Gb)\right)\leq  \mathbb E \phi \left( \frac{\pi}{2} \|a-b\|_2 \langle \nabla f(X), Y\rangle \right),
\end{align*} where $X,Y$ are independent copies of a Gaussian $n$-dimensional random vector. 
\end{theorem} 
The proof is left as an exercise to the interested reader (see also \cite{PV-dvo}).

\smallskip

\noindent 3. For $a,b\in S^{k-1}$ with $\langle a,b\rangle=0$ the above statements are reduced to the inequalities we 
discussed in Section 2.
\end{remarks}

The next result is an application of Theorem \ref{thm:Sch-funct-ineq} for the $\ell_p$ norm.

\begin{theorem} \label{thm:main-5-1}
Let $n$ be large enough and let $2<p <c_0\log n$. Let $a,b \in S^{k-1}$ and let $G=(g_{ij})_{i,j=1}^{n,k}$ be standard Gaussian random variables. Then, 
\begin{align*}
\left( \mathbb E \left| \|Ga\|_p-\|Gb\|_ p\right|^r \right)^{1/r} \lesssim \|a-b\|_2 \psi(n,p,r) \mathbb E\|Z\|_p,
\end{align*} for $r\geq 2$, where $\psi(n,p,r)$ is defined as:
\begin{align*}
\psi(n,p,r):=\sqrt{r} \min \left\{ \frac{1}{\sigma_p n^{1/p}}, \frac{\sigma_{2p-2}^{p-1}}{n^{1/2}\sigma_p^p } \left(1+\frac{pr}{\sigma_{2p-2}^2n^{\frac{1}{p-1}}} \right)^{\frac{p-1}{2}} \right\}
\end{align*}
Moreover, for any $\varepsilon >0$ one has:
\begin{align*}
P \left( \left| \|Ga\|_p-\|Gb\|_ p\right| > \varepsilon \mathbb E\|Z\|_p \right) \leq C \exp \left (-c \tau \left( n,p,\frac{\varepsilon}{\| a-b\|_2} \right) \right),
\end{align*} where 

\begin{align*}
\tau(n,p,t) := \max \left\{ t^2 pn^{2/p} , \min \left\{ \frac{t^2 n}{C^p}, (t n)^{2/p}  \right \} \right \}, \quad t>0
\end{align*}
and $C,c>0$ are absolute constants.
\end{theorem}

\noindent {\it Proof.} In view of Theorem \ref{thm:Sch-funct-ineq} we need an upper estimate for the quantity:
\begin{align} \label{eq:main-5-1-1}
\left( \mathbb E \big\| \nabla \|X\|_p \big\|_2^r \right)^{1/r}
= \left( \mathbb E \frac{\|X\|_{2p-2}^{r(p-1)} }{\|X\|_p^{r(p-1)} } \right)^{1/r} \leq 
\frac{I_{r(p-1)}^{p-1}(\gamma_n,B_{2p-2}^n) }{I_{-r(p-1)}^{p-1}(\gamma_n,B_p^n)},
\end{align} where in the last step we have used Proposition \ref{prop:Harris}. A standard application 
of Lemma \ref{lem:log-sob-moms} (we use \eqref{eq: 2.19}) yields:
\begin{align} \label{eq:main-5-1-2}
\frac{I_{r(p-1)}^{p-1}(\gamma_n,B_{2p-2}^n) }{n^{1/2} \sigma_{2p-2}^{p-1}} 
= \frac{I_{r(p-1)}^{p-1}(\gamma_n,B_{2p-2}^n) }{I_{2p-2}^{p-1}(\gamma_n,B_{2p-2}^n)} 
\leq \left(1+\frac{(p-1)(r-2)}{\sigma_{2p-2}^2n^{\frac{1}{p-1}}} \right)^{\frac{p-1}{2}}.
\end{align} Moreover, from Proposition \ref{prop:small-ball} we see that:
\begin{align} \label{eq:main-5-1-3}
\frac{I_{-r(p-1)}^{p-1} (\gamma_n,B_p^n)}{n^{1-1/p} \sigma_p^{p-1} } 
\gtrsim \frac{I_{-r(p-1)}^{p-1} (\gamma_n,B_p^n)}{I_p^{p-1}(\gamma_n,B_p^n)} \gtrsim 1,
\end{align} for $r \leq c_1k(B_p^n)$.
Plugging estimates \eqref{eq:main-5-1-2} and \eqref{eq:main-5-1-3} in \eqref{eq:main-5-1-1} we find:
\begin{align*}
\left( \mathbb E \big\| \nabla \|X\|_p \big\|_2^r \right)^{1/r} 
\lesssim \frac{ (\sigma_{2p-2} / \sigma_p)^{p-1} }{n^{1/2-1/p} } \left(1+\frac{p(r-2)}{\sigma_{2p-2}^2n^{\frac{1}{p-1}}} \right)^{\frac{p-1}{2}},
\end{align*} for $2\leq r\leq c_1k(B_p^n)$. Taking into account that $\big\| \nabla \|X\|_p \big\|_2 \leq 1$ a.s. we conclude the 
first assertion.
For the distributional inequality we argue as in the proof of Theorem \ref{thm:conc-p>2-I-p}, i.e. we use Markov's inequality and
the previous estimate. \prend

\medskip

\noindent {\bf The chaining method: Dudley-Fernique decomposition.} For each $j=1,2,\ldots$ consider $\delta_j$-nets ${\cal N}_j$ on $S^{k-1}$ with cardinality $|{\cal N}_j| \leq (3/\delta_j)^k$ (see \cite[Lemma 2.6]{MS}). Note that for any
$\theta\in S^{k-1}$ and for  all $j$ there exist $u_j\in {\cal N}_j$ with $\|\theta-u_j\|_2\leq \delta_j$ and by the triangle inequality it follows that
$\|u_j-u_{j-1}\|_2\leq \delta_j+\delta_{j-1}$. Moreover, if we assume that $\delta_j\to 0$ as $j\to \infty$ and $(t_j)$ is sequence of numbers 
with $t_j\geq 0$ and $\sum_j t_j\leq 1$ then, for any $\varepsilon>0$ we have the following:

\smallskip

\noindent {\it Fact.} Set $E:=\mathbb E\|X\|$. If we define the following sets:
\begin{align*}
A:= \left\{\omega \mid \, \exists \theta\in S^{k-1} :  \big| \| G_\omega(\theta)\|-E \big| > \varepsilon E \right\}, \\
 A_1:= \left \{ \omega \mid \exists u_1\in {\cal N}_1 : \big| \|G_\omega(u_1)\| -E \big| > t_1\varepsilon E \right\} \nonumber
 \end{align*} and for $j\geq 2$  
\begin{align*} A_j:= \left \{\omega \mid \exists u_j\in {\cal N}_j, u_{j-1}\in {\cal N}_{j-1} : \left| \|G_\omega(u_j)\|- \|G_\omega(u_{j-1})\| \right| > t_j \varepsilon E \right \},
\end{align*} then one has: $A\subseteq \bigcup_{j=1}^\infty A_j$ (see also \cite{Sch1}).

\medskip

Now we apply the above chaining method for the $\ell_p$ norm with $p>2$ and we employ the distributional inequality 
of Theorem \ref{thm:main-5-1} to prove our second main result:

\begin{theorem} [Random Dvoretzky for $\ell_p^n$]  \label{thm:rdm-dvo-p>2} For all large $n$, for any $1\leq p \leq \infty$ and 
for every $0<\varepsilon<1$ there exists $k(n,p,\varepsilon)$ with the following property: the random $k$-dimensional subspace of 
$\ell_p^n$ with $k\leq k(n,p,\varepsilon)$ is $(1+\varepsilon)$-Euclidean with probability greater than 
$1-C\exp(-c k(n,p,\varepsilon) )$, where $k(n,p,\varepsilon)$ is estimated as follows:
\begin{itemize}
\item [\rm (i)] For $1\leq p<2$ we have:
\begin{align*}
k(n,p,\varepsilon) \gtrsim \varepsilon^2 n,
\end{align*}
\item [\rm (ii)] For $2<p < c_0 \log n$ we have:
\begin{align*}
k(n, p, \varepsilon) \gtrsim \left\{ \begin{array}{lll}
(Cp)^{-p} \varepsilon^2 n, & 0<\varepsilon \leq (Cp)^{p/2} n^{-\frac{p-2}{2(p-1)}} \\
\frac{1}{p} (\varepsilon n)^{2/p} , & (Cp)^{p/2} n^{-\frac{p-2}{2(p-1)}} < \varepsilon \leq 1/p\\
\varepsilon pn^{2/p}/ \log\frac{1}{\varepsilon} , & 1/p < \varepsilon <1.
\end{array} \right.
\end{align*} Moreover, for $p <c_0 \log n$ we have:
\begin{align*}
k(n,p,\varepsilon) \gtrsim \log n/ \log \frac{1}{\varepsilon},
\end{align*}
\item [\rm (iii)] For $c_0 \log n <p \leq \infty$ we have:
\begin{align*}
k(n, p, \varepsilon) \gtrsim \varepsilon \log n /\log \frac{1}{\varepsilon},
\end{align*}
\end{itemize} where $C, c, c_0>0$ are absolute constants.
\end{theorem}

\noindent {\it Sketch of proof.} For $1\leq p<2$ the assertion follows from Theorem \ref{thm: VMil} and the fact that $k(B_p^n) \simeq n$. 
Let $2<p<c_0\log n$ and fix $0<\varepsilon<1/p$. Choose $\delta_j=e^{-j}$, $t_j = s_p^{-1} j^{p/2}e^{-j}$, with $s_p:=\sum_{j=1}^\infty j^{p/2}e^{-j}$. 
Then, according to 
the previous chaining method we may write:
\begin{align*}
P(A) &\leq C|{\cal N}_1| \exp( -c_1 \tau(n,p,\varepsilon t_1) ) + C\sum_{j=2}^\infty |{\cal N}_{j-1}| \cdot |{\cal N}_j| \exp (- c_1\tau(n,p, \varepsilon  s_p^{-1} t_je^j/4)) \\
&\leq C\sum_{j=1}^\infty (3 e^{j})^{2k} \exp(-c_2 \tau(n,p, s_p^{-1} \varepsilon j^{p/2}) ),
\end{align*} where $\tau(n,p,t)$ was defined in Theorem \ref{thm:main-5-1}, hence:
\begin{align*}
\tau(n,p, t) \simeq \min \left\{ \frac{t^2n}{C_1^p}, (tn )^{2/p}\right\}, \, t>0.
\end{align*} Note that
\begin{align*}
\tau(n, p, s_p^{-1} \varepsilon j^{p/2}) \gtrsim j \min\left\{ \frac{\varepsilon^2 n}{(Cp)^p}, \frac{(\varepsilon n)^{2/p}}{p} \right\}= j k(n, p, \varepsilon),
\end{align*} where we have used the fact that $s_p \lesssim \sqrt{p} (\frac{p}{2e})^{p/2}$. Therefore, we have:
\begin{align*}
P(A) &\leq C \sum_{j=1}^\infty \exp \left( c_3j k- c_4 j k(n,p,\varepsilon) \right) \\
&\leq \sum_{j=1}^\infty \exp \left( - \frac{c_4}{2} j k(n,p,\varepsilon) \right) \leq C' \exp \left( - \frac{c_4}{2} k(n,p,\varepsilon) \right) .
\end{align*}
as long as $k\leq \frac{c_4}{2c_2} k(n,p, \varepsilon)$. 

In the case that $p<c_0 \log n$ and $p\gg 1$ for the range $1/p<\varepsilon<1$ we have for any fixed $\theta\in S^{k-1}$ 
the concentration inequality
\begin{align*}
P\left( \big| \|G\theta\|_p -\mathbb E\|X\|_p \big| > \varepsilon \mathbb E\|X\|_p \right) \leq C \exp(-c\varepsilon k(B_p^n) ),
\end{align*} by Proposition \ref{prop: weak-conc}. Thus, the classical net argument yields the estimate: 
$k(n,p,\varepsilon) \gtrsim \varepsilon k(B_p^n) /\log \frac{1}{\varepsilon}$. We omit the details.

Moreover, for $2< p < c_0 \log n$ but $p\simeq \log n$, the main result of 
Section 2 shows that ${\rm Var}\|X\|_p \lesssim n^{-c_1}$ for some absolute constant $c_1>0$. Therefore, Chebyshev's probabilistic 
inequality and the net argument as before implies $k(n, p, \varepsilon) \gtrsim \log n / \log \frac{1}{\varepsilon}$. 

Finally, for $p\gtrsim \log n$ we employ Corollary \ref{cor:conc-large-p} combined with the net argument again to 
get $k(n,p,\varepsilon) \gtrsim \varepsilon \log n/ \log \frac{1}{\varepsilon}$.  \prend

\medskip

Below we show that the dependence on $\varepsilon$ we get for the randomized Dvoretzky in $\ell_p^n$, for fixed $2< p<\infty$ is 
essentially optimal. We have the following:

\begin{theorem} [Optimality in the Random Dvoretzky for $\ell_p^n$] \label{thm:sharp-dvo}
Let $2< p< c_0 \log n$. Assuming that with probability larger than $1-e^{-\beta k}$, a $k$-dimensional subspace satisfies
that the ratio between the $\ell_p^n$ norm and a multiple of the $\ell_2^n$ norm is $(1+\varepsilon)$ equivalent for all 
vectors in the subspace, with $\frac{2^{p/2}}{p} n^{-\frac{p-2}{2(p-1)}} < \varepsilon <1/p$, then 
$k\lesssim  \beta^{-1} \varepsilon^{2/p} k(B_p^n)$.
\end{theorem}

For the proof we will need the next lemma from \cite{Sch2}:

\begin{lemma} \label{lem:meas-grass-gauss}
Let $1\leq k\leq n-1$ and let ${\cal A}\subset G_{n,k}$ be a $\nu_{n,k}$-measurable set. Then, for $U_{\cal A}: =\bigcup \{F \mid F\in {\cal A} \}$ we
have: 
\begin{align*}
\nu_{n,k}({\cal A}) \leq [\gamma_n(U_{\cal A})]^k.
\end{align*}
\end{lemma}

\noindent {\it Proof of Theorem \ref{thm:sharp-dvo}.} 
Let $0<\varepsilon<1/3$ and define the 
collection of all $k$-dimensional subspaces of a space $(\mathbb R^n, \| \cdot \|)$ for which the restricted norm there
has distortion (with respect to the Euclidean norm) at most $1+\varepsilon$:
\begin{align*}
{\cal A}_\varepsilon:=\{F\in G_{n,k} \mid \exists \lambda_F : \lambda_F \leq \|\theta\| \leq (1+\varepsilon) \lambda_F, \, \forall \theta\in S_F\}.
\end{align*} Note that for $F\in {\cal A}_\varepsilon$ we have: $(1+\varepsilon)^{-1}M_F \leq \lambda_F\leq M_F$. Thus instead of working with 
$\lambda_F$ we may define ${\cal A}_\varepsilon$ using $M_F:=M(F\cap B)$ (here $B=\{x : \|x\|\leq 1\}$) namely, if
\begin{align*}
{\cal F}_\varepsilon:= \left\{ F\in G_{n,k} \mid (1+\varepsilon)^{-1} M_F\leq \|\theta\| \leq (1+\varepsilon)M_F\, \; \forall \theta\in S_F \right\},
\end{align*} then we get ${\cal A}_\varepsilon \subset {\cal F}_{\varepsilon}$. Define further:
\begin{align*}
{\cal B}_\varepsilon:= \left\{ F\in {\cal F}_\varepsilon \mid (1-2\varepsilon)\frac{\mathbb E\|g\|}{\mathbb E\|g\|_2}\leq 
M_F\leq (1+2\varepsilon)\frac{\mathbb E\|g\|}{\mathbb E\|g\|_2}\right\}
\end{align*} and note that ${\cal F}_\varepsilon, \; {\cal B}_\varepsilon$ are measurable.\footnote{The map $F\mapsto M_F$ is Lipschitz continuous 
with respect to the unitarily invariant metric $d$ on $G_{n,k}$ defined as: $d(E,F)= \inf\{ \|I-U\|_{\rm op} : U(E)=F, \, U\in O(n)\}, \; E,F\in G_{n,k}$.} 
Hence, an application of Lemma \ref{lem:meas-grass-gauss} yields:
\begin{align*}
\nu_{n,k}({\cal F}_\varepsilon) &=  \nu_{n,k}({\cal F}_\varepsilon \setminus {\cal B}_\varepsilon)+ \nu_{n,k}({\cal B}_\varepsilon) \\
& \leq \left[ \gamma_n \left( \left\{ x: \|x\| \geq \frac{1+ 2 \varepsilon}{1+\varepsilon}\frac{\mathbb E\|g\|}{\mathbb E\|g\|_2} \|x\|_2  \; {\rm \bf or} \; 
\|x\| \leq (1+ \varepsilon)(1- 2\varepsilon)\frac{\mathbb E\|g\|}{\mathbb E\|g\|_2} \|x\|_2 \right\} \right) \right]^k + \\
& \left[ \gamma_n \left( \left\{ x: \frac{1-2\varepsilon}{1+\varepsilon} \|x\|_2 \frac{\mathbb E\|g\|}{\mathbb E\|g\|_2} \leq \|x\|\leq (1+\varepsilon)(1+2\varepsilon)  \frac{\mathbb E \|g\|}{\mathbb E \|g\|_2} \|x\|_2 \right\}\right) \right]^k.
\end{align*} Apply this argument for the $\ell_p$ norm with $2<p< c_0 \log n$ and consider the next claim which follows easily by
Theorem \ref{thm:conc-1<p<2} and Proposition \ref{prop:conc-p>2-u-l-b}: 

\smallskip

\noindent {\it Claim.} For every $2^{p/2} p^{-1} n^{-\frac{p-2}{2(p-1)} } < t < 1/p$ we have:
\begin{align*}
c e^{-Cp (tn)^{2/p} } \leq P \left( \|g\|_p \leq \frac{(1-t)\mathbb E \|g\|_p}{\mathbb E\|g\|_2} \|g\|_2 \; {\rm \bf or} \;  
\|g\|_p \geq \frac{(1+t)\mathbb E \|g\|_p}{\mathbb E\|g\|_2} \|g\|_2 \right) \leq C e^{-c (tn)^{2/p} }.
\end{align*}

\smallskip

\noindent Now assume that $2^{p/2} p^{-1} n^{-\frac{p-2}{2(p-1)} }<\varepsilon<1/p$, so by the previous claim we get:
\begin{align*}
\nu_{n,k}({\cal F}_\varepsilon) &\leq C^k e^{-ck(\varepsilon n)^{2/p}}+ (1-ce^{-Cp(\varepsilon n)^{2/p}})^k 
\leq e^{-c'k(\varepsilon n)^{2/p}} +1-ce^{-Cp(\varepsilon n)^{2/p}}.
\end{align*}
Now employing the assumption that $\nu_{n,k}({\cal F}_\varepsilon) \geq 1- e^{-\beta k}$ for some absolute constant $\beta>0$ and
that $\beta \ll ( \varepsilon n)^{2/p}$, we conclude:
\begin{align*} 
1-ce^{-Cp(\varepsilon n)^{2/p}} \geq 1-e^{-\beta k} -e^{-c'k(\varepsilon n)^{2/p}}  \geq 1- 2e^{-c'' \beta k},
\end{align*} which implies $k\leq \frac{C'}{\beta} p (\varepsilon n)^{2/p}$, as required. \prend

\section{Further remarks and questions}

\noindent {\bf 1. Instability of the variance.} It is worth mentioning that the variance is not an isomorphic invariant.  
One can observe that:
\begin{quote} \it There exists absolute constant $0< c_0 <1$ with the following property: for every $n\geq 2$
there exist 1-symmetric convex bodies $K$ and  $L$ on $\mathbb R^n$ such that:
\begin{align*}
{\rm Var}\|Z\|_K \simeq \frac{1}{n^{\delta} \log n}, \quad {\rm Var}\|Z\|_L \simeq \frac{1}{\log n} \quad {\rm and} 
\quad e^{-1/c_0}L \subseteq K \subseteq L,
\end{align*} where $\delta =1-c_0 \log 2$ and $Z\sim N({\bf 0},I_n)$.
 \end{quote}  
 
\noindent Indeed; for $p_0:= c_0 \log n$, where 
$0<c_0<1$ as in Theorem \ref{thm:var-ell_p}, we consider the bodies $K:=B_{p_0}^n$ and $L:=B_\infty^n$. We can easily see that
these bodies enjoy the aforementioned properties.

\medskip

\noindent {\bf 2. Non-centered moments.} We know that for any centrally symmetric convex body $T$ on $\mathbb R^n$ one has: 
\begin{align*}
\frac{c_1r}{n} \leq \left( \frac{I_r(\gamma_n,T)}{I_1(\gamma_n,T)} \right)^2 -1 \leq \frac{c_2r}{k(T)},
\end{align*} for all $r\geq 2$, where $c_1,c_2>0$ are absolute constants. This follows from the lower estimate in \eqref{eq:2.15} and 
Lemma \ref{lem:log-sob-moms}.  In particular, for $1\leq r \leq k(T)$ we obtain:
\begin{align} \label{eq:6.2}
\frac{c_1' r}{n} \leq \frac{I_r(\gamma_n,T)}{I_1(\gamma_n,T)}  -1 \leq \frac{c_2' r}{k(T)}
\end{align} and
when $k(T) \simeq n$ we readily see that this estimate is sharp up to constants, in particular for the $\ell_p$ norms with $1\leq p\leq 2$. Furthermore, one can show that the same behavior holds true for $2< p< c_0\log n$ even though the critical dimension $k(B_p^n)$ 
in that case is much smaller than $n$. For 
$2<p<c_0 \log n$ we have:
\begin{align*}
\frac{I_r(\gamma_n, B_p^n)}{I_1(\gamma_n,B_p^n)} \leq 1 +\frac{C^p}{n}r,
\end{align*} for all $1\leq r\leq k(B_p^n)/C$. In fact for the negative moments this is already clear if we take into account Theorem 
\ref{thm:stability-moms}, Proposition \ref{prop: stab-ell_p} and Theorem \ref{thm:stability-r-means}. More precisely 
we have:
For $1\leq p < c_0\log n$ and for any $1\leq r\leq c k(B_p^n)$ we get:
\begin{align*}
\max \left \{ \frac{I_1(\gamma_n, B_p^n)}{I_{-r}(\gamma_n,B_p^n)} , \, \frac{I_r(\gamma_n, B_p^n)}{I_1(\gamma_n,B_p^n)} \right\} 
\leq 1 +\frac{C^p}{n}r
\end{align*} and for $p\geq c_0\log n$ and $1\leq r \leq c k(B_p^n)$ we have:
\begin{align*}
\max \left \{ \frac{I_1(\gamma_n, B_p^n)}{I_{-r} (\gamma_n,B_p^n)}, \, \frac{I_r(\gamma_n, B_p^n)}{I_1(\gamma_n,B_p^n)}  \right\}
\leq 1 +\frac{C}{(\log n)^2}r.
\end{align*} We should note here the next threshold phenomenon when $2<p\leq \infty$: 

\begin{itemize}
\item $2<p\leq c_0\log n$: It is $I_r/I_1 -1 \lesssim_p r / n  =O_p(n^{2/p-1})$ for $1 \leq r\leq c_1 k(B_p^n)$ while for $r\geq c_2 k(B_p^n)$ we have $I_r/I_1 -1 \simeq 1$.
\item $p>c_0\log n$: It is $I_r/I_1 -1 \lesssim r / (\log n)^2 =O( (\log n)^{-1})$ for $1\leq r\leq c_1 k(B_p^n)$, while for $r\geq c_2k(B_p^n)$ we have $I_r/I_1 -1 \simeq 1$,
\end{itemize} for absolute constants $0< c_1< c_2$. The detailed study of this phenomenon will be presented elsewhere. Let us also note
that although the behavior of the quantities $\frac{I_r}{I_1} -1, \, \frac{I_1}{I_{-r}}-1$ is completely determined for 
the $\ell_p$ norms -- it is of the order $r/n$ for $1\leq r\leq ck(B_p^n)$ -- combining this information with Markov's inequality we still do not
derive the optimal concentration inequality in the whole range $2<p <\infty$.

\medskip

\noindent {\bf 3. Gaussian concentration and randomized Dvoretzky.} One can show that the Gaussian concentration for norms 
$\|\cdot\|_A$ with $k(A) \simeq n$ is essentially optimal:

\begin{lemma} Let $\alpha\in (0,1)$ and let $A$ be centrally symmetric convex body on 
$\mathbb R^n$ with $k=k(A) \geq \alpha n $. Then, 
\begin{align*}
P \left( \big | \|Z\|_A - \mathbb E \|Z\|_A \big |\geq \varepsilon \mathbb E\|Z\|_A \right) \geq c e^{-C\varepsilon^2 k / \alpha^2},
\end{align*} for all $n^{-1/2} <\varepsilon<1$.
\end{lemma}

\noindent {\it Proof.} Set $I_q^q=\mathbb E \|Z\|_A^q$. Taking into account \eqref{eq:6.2} we may write:
\begin{align*}
1+ \frac{c_1r}{n} \leq \frac{I_r}{I_1} \leq \sqrt{ 1+\frac{C_1r}{k} },
\end{align*} for all $r\geq 2$. Let $n^{-1/2} <\varepsilon < 1$. If we set $r_0:=\frac{2n\varepsilon}{c_1}$, then by previous estimates
and the Paley-Zygmund inequality we have:
\begin{align*}
P\left( \|Z\| > (1+\varepsilon) I_1 \right)\geq P\left( \|Z\| > \frac{1+\varepsilon}{1+\frac{c_1r_0}{n}} I_{r_0}\right) = 
P( \|Z\|> \delta I_{r_0}) \geq (1-\delta^{r_0})^2 \frac{I_{r_0}^{2r_0}}{I_{2r_0}^{2r_0}} \geq c_2 e^{-C_2 r_0^2/k } ,
\end{align*} where $\delta := \frac{1+\varepsilon}{1+\frac{c_1r_0}{n}}$. The result easily follows. \prend

\smallskip

Although the Gaussian concentration for spaces $E=(\mathbb R^n, \|\cdot\|)$ with $k(E) \simeq n$ is sharp, the argument provided in Section 5 fails to give the optimal dependence on $\varepsilon$ in randomized Dvoretzky. 
The reason is that in Gauss' space, norms with concentration estimate less than $e^{-\varepsilon^2n}$ cannot be distinguished 
from the Euclidean norm. Therefore it is more appropriate to work on the sphere when we study almost spherical sections in normed spaces. 

\medskip

\noindent {\bf 4. Refined Gaussian concentration and "new dimensions".} The reader should notice that the refined form of 
the Gaussian concentration for $2<p<\infty$ (Theorem \ref{thm:conc-full}) and moreover Theorem \ref{thm:rdm-dvo-p>2} provide random, 
almost Euclidean 
subspaces of relatively large dimensions in which the norm has very small distortion. Previously, that phenomenon could not be observed 
if one was using the classical concentration inequality in terms 
of the Lipschitz constant. In order to illustrate
this let us consider an example, say the $\ell_p$ norm with $p=5$. The classical setting yields the 
existence of random $k$-dimensional sections
of $B_5^n$ which are $(1+\varepsilon)$-isomorphic to a multiple of $B_2^k$ as long as $k \lesssim \varepsilon^2 n^{2/5}$. The latter 
is relatively large when $\varepsilon \gg n^{-1/5}$. Now, we may consider distortions smaller than $n^{-1/5}$, in fact as small as $n^{-1/2}$,
since $\tau(n,5,\varepsilon) \simeq \min\{\varepsilon^2 n, (\varepsilon n)^{2/5} \}$. For instance (for $\varepsilon \simeq n^{-2/5}$) 
the random $k$-dimensional section of $B_5^n$ with $k\simeq n^{1/5}$, is $(1+n^{-2/5})$-isomorphic to a multiple of $B_2^k$ with probability greater than $1-e^{-cn^{1/5}}$.

\medskip

\medskip

\noindent {\bf 5. The existence of $\log(1/\varepsilon)$ as $p\to \infty$.} Note that Theorem \ref{thm:stability-r-means} and 
furthermore
Corollary \ref{cor:conc-large-p} suggest that the concentration of the $\ell_p$ norm with $p\gtrsim \log n$ is similar 
with the one we get for the $\ell_\infty$ norm. This
means that the classical net argument yields random subspaces which are $(1+\varepsilon)$-spherical as 
long as $k\lesssim \varepsilon \log n/ \log \frac{1}{\varepsilon}$. We do not know if this $\log(1/\varepsilon)$ term is needed, 
for this range of $p$. As an easy corollary of the main result of \cite{Tik}, we have:

\begin{proposition}
Let $p> (\log n)^2$ and $\varepsilon \in(0,1/3)$. If the random $k$-dimensional subspace of $\ell_p^n$ is 
$(1+\varepsilon)$-spherical with probability greater than $3/4$, then $k\leq C \varepsilon \log n/ \log \frac{1}{\varepsilon}$, where
$C>0$ is an absolute constant.
\end{proposition}

\medskip

\noindent {\bf Acknowledgements.} The authors are indebted to the anonymous referees whose valuable comments
helped to improve the presentation of this note. The authors would also like to thank Konstantin Tikhomirov who was interested 
in their question and kindly allowed them to include his argument here.

\appendix

\section{An anti-concentration estimate by Tikhomirov}

After the second named author presented the results of the paper, to the Functional Analysis Seminar 
at the Math Department in University of Alberta, Tikhomirov was motivated to study the lower bound 
of the variance for $p \gtrsim \log n$. He proved that the upper bound $(\log n)^{-1}$ is tight. 
Moreover, he proved that the concentration we obtain in Corollary \ref{cor:conc-large-p} is sharp. 
We are indebted to him for kindly allowing us 
to present his argument here:

\begin{theorem} [Tikhomirov, 2016] Let $p \geq C_0 \log n$. Then, one has:
\begin{align}
P\left( \big| \|Z\|_p -\mathbb E \|Z\|_p \big| >\varepsilon \mathbb E \|Z\|_p \right) \geq ce^{-C\varepsilon \log n},
\end{align} for all $\varepsilon \in (0,1)$, where $c,C, C_0>0$ are absolute constants and $Z \sim N({\bf 0},I_n)$.
\end{theorem}

In fact something more is true. In order to formulate it we need a little bit of terminology. Let $(\Omega,\Sigma, P)$ be the probability space.
In what follows, we let $X=(x_1,x_2,\dots,x_n)=(|g_1|,|g_2|,\dots,|g_n|)$, where $g_1,g_2,\dots, g_n$ are i.i.d. standard
Gaussian variables. Further, for a random variable $\eta$ let $\mathcal Q (\eta,t)$ be its L\'{e}vy concentration function defined by
$$\mathcal Q(\eta,t):=\sup_{\lambda \in \mathbb R} P \left ( |\eta - \lambda | \leq t \right ), \; t>0.$$
By $z^*$ we denote the non-increasing rearrangement of a vector $z$.

\begin{proposition}\label{main-A}
There is a universal constant $C>0$ such that for any $n\geq C$, $p\geq 12\log n$
and any $\varepsilon\in(0,1)$ we have
$$\mathcal Q \left( \|X\|_p, \varepsilon\sqrt{\log n} \right)\leq 1-0.07n^{-120\varepsilon}.$$
\end{proposition}

By $o(1)$ we mean any quantity which is a function of $n$ and
becomes arbitrarily small for large enough $n$.
The dimension $n$ is always assumed to be large.
Further, for any $s\in (0,1)$ let $\xi_s$ be the quantile of order $s$ with respect to the distribution of $|g_1|$, i.e.
the number satisfying the equation
$$\sqrt{\frac{2}{\pi}} \int_0^{\xi_s} \exp(-t^2/2) \, dt = P \left( |g_1| \leq \xi_s \right) = s.$$
By $y=(y_1,y_2, \dots, y_n)$ we denote a non-random vector of quantiles, where
$$y_i:=\xi_{1-(i-0.5)/n}, \; i=1,2,\dots,n.$$

\noindent  It follows from Lemma \ref{lem:Gordon-ineq} that
\begin{equation}\label{basic est}
\left( a^{-1}-a^{-3} \right) \exp(-a^2/2) < \int_a^\infty \exp(-t^2/2) \, dt < a^{-1} \exp(-a^2/2), \; a>0.
\end{equation} 
Using this estimate, it is elementary to check that
$$y_1=\left( 1+o(1) \right) \sqrt{2\log n}.$$

\begin{lemma}\label{l quantile int}
We have $y_1^2-y_i^2\geq\log i$ for all $1\leq i\leq 0.317n$.
\end{lemma}

\noindent {\it Proof.} By the definition of the quantiles, we have
$$\int_{y_m}^\infty \exp(-t^2/2) \, dt = \sqrt{\frac{\pi}{2}}\frac{2m-1}{2n},  \; m=1,2,\dots,n.$$
Together with \eqref{basic est}, it gives:
$$\left( y_1^{-1}-y_2^{-3} \right) \exp(-y_1^2/2) < \frac{y_i^{-1}}{2i-1} \exp(-y_i^2/2),$$
whence
$$\frac{y_1}{y_i} \exp \left( y_1^2/2-y_i^2/2 \right) > (1-o(1))(2i-1).$$
It can be checked that, under the assumption that $y_i\geq 1$ (which holds true since $i\leq 0.317 n$), we have
$$\frac{y_1}{y_i} \exp \left( y_1^2/2-y_i^2/2 \right) \leq \exp \left( y_1^2-y_i^2 \right).$$
Plugging the estimate into the previous formula and using the rough bound $(1-o(1))(2i-1)\geq i$ for $i\geq 2$, we obtain the result.
\prend

\medskip

The next lemma is checked by a direct computation:
\begin{lemma}\label{l bound x1star}
Denote $y_0:=\xi_{1-1/(4n)}$. Then
$$P \left( x_1^*\in [y_1,y_0] \right) \geq 0.17.$$
\end{lemma}

\begin{lemma}\label{l bound xistar}
For every $i\geq e^2$ we have
$$P \left( x_i^*\geq y_{[i/e^2]} \right) \leq \exp(-i).$$
\end{lemma}

\noindent {\it Proof.} Recall that by Chernoff's inequality we have for any $i$ and any $s\in (1-i/n,1)$:
\begin{align*}
\sum_{j=0}^{n-i} {n\choose j}(1-s)^{n-j}s^j
\leq \exp \left( (n-i) \log \frac{sn}{n-i}+ i \log \frac{(1-s)n}{i} \right) \leq \exp(i) \left( \frac{(1-s)n}{i} \right)^i.
\end{align*}
Hence, denoting $s:=1-\frac{[i/e^2]-1/2}{n}$, we get
\begin{align*}
P\left( x_i^*\geq y_{[i/e^2]} \right) = P( x_i^*\geq \xi_s) =\sum_{j=0}^{n-i} {n\choose j} P ( |g_1| \geq \xi_s)^{n-j} P( |g_1|\leq \xi_s)^j 
&=\sum_{j=0}^{n-i} {n\choose j} (1-s)^{n-j} s^j \\
&\leq \exp(i) \left( \frac{(1-s)n}{i} \right)^i \\
&\leq \exp(-i).
\end{align*}
\prend

\noindent  Let us fix any $\varepsilon\in (0,1]$ and denote
\begin{align*}
Q_1 &:= \left \{ (z_1,\dots,z_n) \in \mathbb R_+^n : \, z_1^* \in [y_1,y_0] ; \; z_i^* \leq y_{[i/e^2]} \; \forall \; i \geq 2 \right \}; \\
Q_2 &:= \left \{ (z_1,\dots,z_n) \in \mathbb R_+^n : \, z_1^* \in [y_1+60\varepsilon \sqrt{\log n}, y_0+60\varepsilon \sqrt{\log n}]
; \; z_i^* \leq y_{ [i/e^2]} \; \forall \; i \geq 2 \right \}.
\end{align*}

Further, let $\mathcal E_k:= \left\{ \omega \in \Omega : \, X(\omega) \in Q_k \right \}$ ($k=1,2$).
It is easy to see from Lemmas~\ref{l bound x1star} and~\ref{l bound xistar} that
$P( \mathcal E_1) \geq 0.15$ (note that for $i<e^2$ the condition $x_i^*\leq y_{[i/e^2]}=y_0$ is fulfilled automatically
provided that $x_1^*\leq y_0$).

\begin{lemma}\label{l pnorm}
Let $p\geq 12\log n$ and $z\in Q_1$. Then, 
$$\|z\|_p^p\leq 3e^2(z_1^*)^p.$$
\end{lemma}

\noindent {\it Proof.} We have
$$\sum_{i=2}^n (z_i^*)^p
\leq (e^2-1)(z_1^*)^p+\sum_{k=2}^{[n/e^2]+1} \sum_{i=[e^2(k-1)]+1}^{[e^2 k]}(z_i^*)^p
\leq (e^2-1)(z_1^*)^p+(e^2+1) \sum_{k=2}^{[n/e^2]+1}y_{k-1}^p.$$
Applying Lemma~\ref{l quantile int}, we get
$$\sum_{k=2}^{[n/e^2]}y_{k-1}^p\leq \sum_{k=2}^{[n/e^2]+1}(y_1^2-\log (k-1))^{p/2}
\leq y_1^p\sum_{k=2}^{\infty}(k-1)^{-p/(4\log n)}<1.21y_1^p.$$
Finally, we get
$$\|z\|_p^p\leq e^2(z_1^*)^p+(e^2+1)y_1^p\leq 3e^2(z_1^*)^p.$$
\prend

Next, consider an operator $T:\mathbb R_+^n \to \mathbb R_+^n$ which acts by
adding $60\varepsilon\sqrt{\log n}$ to the largest coordinate of a vector, i.e.\
for $z=(z_1,z_2,\dots,z_n)\in \mathbb R_+^n$ with $k= \min \left \{ s : s={\rm argmax} \{z_1,z_2,\dots,z_n\} \right\}$ we have
$Tz=(z_1,z_2,\dots,z_k+60\varepsilon\sqrt{\log n},\dots,z_n)$.
Obviously, $T$ maps the set $Q_1$ into $Q_2$.

\begin{lemma}\label{l distance}
Let $z\in Q_1$ and $p\geq 12\log n$. Then, 
$$\|Tz\|_p-\|z\|_p > 2\varepsilon\sqrt{\log n}.$$
\end{lemma}

\noindent {\it Proof.} By Lemma~\ref{l pnorm}, we have
$$\sum_{i=2}^n(z_i^*)^p\leq (3e^2-1)(z_1^*)^p.$$
On the other hand, $\left( (Tz)_1^* \right)^p= \left ( 1+60\varepsilon \sqrt{\log n}/z_1^* \right)^p(z_1^*)^p$. Thus,
$$\frac{\|Tz\|_p^p}{\|z\|_p^p}\geq 1+\frac{ \left( 1+60\varepsilon\sqrt{\log n}/z_1^* \right)^p-1}{3e^2}
\geq \left( 1+60 \varepsilon \sqrt{\log n}/z_1^* \right)^{p/(3e^2)},$$
whence
$$\frac{\|Tz\|_p}{\|z\|_p}\geq \left( 1+60 \varepsilon\sqrt{\log n}/z_1^* \right)^{1/(3e^2)}
\geq 1+\frac{2 \varepsilon \sqrt{\log n} }{z_1^*}
> 1+\frac{2 \varepsilon \sqrt{\log n} }{ \|z\|_p}.$$
\prend

\noindent {\it Proof of Proposition~\ref{main-A}.} Denote by $\rho(z)$ ($z\in \mathbb R_+^n$) the probability density function of the vector $X$.
Then from the definition of $T$ and $Q_1$ we have for any $z\in Q_1$:
$$\rho(z)\geq\rho(Tz)\geq \frac{\exp \left( - (y_0+60\varepsilon \sqrt{\log n} )^2/2 \right) }{ \exp(-y_0^2/2)} \rho(z)\geq
n^{-120\varepsilon}\rho(z).$$
Fix any $\lambda\in \mathbb R$. Then from Lemma~\ref{l distance} it follows that for any $z\in Q_1$ we have
$$\max \left\{ \big| \|z\|_p- \lambda \big| , \big| \|Tz\|_p-\lambda \big| \right\} > \varepsilon \sqrt{\log n}.$$
Denote $W_\lambda:= \left \{ z \in \mathbb R_+^n : \, \big| \|z\|_p-\lambda \big| > \varepsilon \sqrt{\log n} \right \}$ and
$\widetilde W_\lambda := \left \{ z \in \mathbb R_+^n : \, \big| \|Tz\|_p-\lambda \big| > \varepsilon \sqrt{\log n} \right \}$. Note that
$T(\widetilde W_\lambda)\subset W_\lambda$ and $T$ is volume preserving transformation, therefore
\begin{align*}
\int_{W_\lambda} \rho(z) \, dz \geq \int_{T(\widetilde W_\lambda)} \rho(z)\, dz = \int_{\widetilde W_\lambda} \rho(Tz) \, dz.
\end{align*}
Moreover, from Lemma~\ref{l distance} it follows that $Q_1\subset W_\lambda\cup\widetilde W_\lambda$, hence we may write:
\begin{align*}
P \left( \big| \|X\|_p-\lambda \big|> \varepsilon \sqrt{\log n} \right) =\int_{W_\lambda} \rho(z) \, dz 
&\geq \frac{1}{2} \left[ \int_{W_\lambda} \rho(z) \, dz + \int_{\widetilde W_\lambda} \rho(T z) \, dz \right] \\
&\geq \frac{1}{2} \int_{W_\lambda\cup\widetilde W_\lambda} \min \left\{ \rho(z), \rho(Tz) \right\} \, dz \\
&\geq \frac{1}{2} \int_{Q_1} \min \left\{ \rho(z), \rho(Tz) \right \} \, dz \\
&\geq \frac{1}{2} n^{-120\varepsilon} \int_{Q_1} \rho(z) \, dz \\
&= \frac{1}{2} n^{-120\varepsilon} P(\mathcal E_1) \geq 0.07n^{-120\varepsilon}.
\end{align*}
\prend


\medskip

\vspace{.5cm} \noindent 

\begin{minipage}[l]{\linewidth}
  Grigoris Paouris: {\tt grigoris@math.tamu.edu}\\
  Department of Mathematics, Mailstop 3368\\
  Texas A\&M University\\
  College Station, TX 77843-3368\\
  
  \medskip
  
  Petros Valettas: {\tt valettasp@missouri.edu}\\
  Mathematics Department\\
  University of Missouri\\ 
  Columbia, MO 65211\\
  
  \medskip
  
  Joel Zinn: {\tt jzinn@math.tamu.edu}\\
  Department of Mathematics, Mailstop 3368\\
  Texas A\&M University\\
  College Station, TX 77843-3368\\

\end{minipage}


\end{document}